\documentclass[11pt]{article}
\usepackage[utf8]{inputenc}
\usepackage[T1]{fontenc}
\usepackage{amsfonts}
\usepackage{amsfonts}
\usepackage{amsfonts}
\usepackage{mathrsfs}
\usepackage{amssymb,amsmath} 
\usepackage[colorlinks,
            linkcolor=red,
            anchorcolor=green,
            citecolor=blue
            ]{hyperref}

 \allowdisplaybreaks

\catcode`!=11
\let\!int\int \def\int{\displaystyle\!int}
\let\!lim\lim \def\lim{\displaystyle\!lim}
\let\!sum\sum \def\sum{\displaystyle\!sum}
\let\!sup\sup \def\sup{\displaystyle\!sup}
\let\!inf\inf \def\inf{\displaystyle\!inf}
\let\!cap\cap \def\cap{\displaystyle\!cap}
\let\!max\max \def\max{\displaystyle\!max}
\let\!min\min \def\min{\displaystyle\!min}
\let\!frac\frac \def\frac{\displaystyle\!frac}
\catcode`!=12

\let\oldsection\section
\renewcommand\section{\setcounter{equation}{0}\oldsection}

\allowdisplaybreaks

\newtheorem{defn}{Definition}[section]
\newtheorem{thm}{Theorem}[section]

\begin{document}
\title{\bf Global well-posedness of a nonlinear Boussinesq-fluid-structure interaction system with large initial data}
\author {Jie Zhang$^{1}$ \thanks{ E-mail:zhangjiee@qztc.edu.cn },\quad Shu Wang$^{2}$ \thanks{Corresponding author. E-mail: wangshu@bjut.edu.cn},\quad Lin Shen$^{3}$ \thanks{ E-mail:shenlin@huanghuai.edu.cn }\\
\small  $^{1}$School of Mathematics and Computer Science, Quanzhou Normal University, \\Quanzhou 362001, P.R. China\\
\small  $^{2}$School of Mathematics, Statistics and Mechanics, Beijing University of Technology,\\ Beijing 100124, P.R. China\\
\small  $^{3}$Department of Mathematics, Huanghuai University,\\ Zhumadian 463000, P. R. China}

\date{}
\maketitle \noindent{\bf Abstract.} In this paper, we consider the global well-posedness of the initial-boundary value problem to a nonlinear Boussinesq-fluid-structure interaction system, which describes the motion of an incompressible Boussinesq-fluid surrounded by an elastic structure with the heat exchange and is one coupled incompressible Boussinesq equations with the wave equation and heat equation by physical interface boundary conditions. Firstly, the global existence of weak solutions to this problem in two/three-dimension is proven by introducing one class of its suitable weak solution and using the compactness method. Then, the uniqueness of the weak solution to this problem in two-dimension is established. Finally, the existence and uniqueness of the global strong and smooth solution to this problem in two-dimension is obtained for any smooth large initial data under the assumptions of suitable compatibility conditions by establishing a priori higher order derivative estimates.
\\[2mm]
\\[2mm]


\section{Introduction}
\quad\quad We study in this paper the global-well-posedness theory of the following initial-boundary value problem to a nonlinear Boussinesq-fluid-structure interaction system
\begin{align}\label{1.1}
\left \{
\begin{array}{lr}
   u_t+ u\cdot\nabla u+\nabla p=\varepsilon\Delta u+\rho e+f_{1},&(x,t)\in\Omega_{f}\times(0,T),\\
  \text{div} u=0,&(x,t)\in\Omega_{f}\times(0,T),\\
  \rho_{t}+u\cdot\nabla\rho=k_{1}\Delta\rho+f_{2},&(x,t)\in\Omega_{f}\times(0,T),\\
  w_{tt}=\mu\Delta w+f_{3},&(x,t)\in\Omega_{s}\times(0,T),\\
  \theta_{t}=k_{2}\Delta \theta+f_{4},&(x,t)\in\Omega_{s}\times(0,T),\\
  w_{t}=u,&(x,t)\in\Gamma\times(0,T),\\
  \varepsilon\frac{\partial u}{\partial n}-pn-\frac{1}{2}(u\cdot n)u=\mu\frac{\partial w}{\partial n},&(x,t)\in\Gamma\times(0,T),\\
  \rho=\theta,&(x,t)\in\Gamma\times(0,T),\\
  k_{1}\frac{\partial \rho}{\partial n}-\frac{1}{2} (u\cdot n)\rho=k_{2}\frac{\partial \theta}{\partial n},&(x,t)\in\Gamma\times(0,T),\\
  w=0,\frac{\partial \theta}{\partial N}=0, &(x,t)\in\Gamma_{out}\times(0,T),\\
  (u,\rho)(x,0)=(u_{0},\rho_{0})(x),&x\in\Omega_{f},\\
  (w,w_{t},\theta)(x,0)=(w_{0},w_{1},\theta_{0})(x),&x\in\Omega_{s},
\end{array}
\right.
\end{align}
where the unknown $u, p, \rho, w$ and $\theta$ denote the fluid velocity, the fluid pressure, the fluid scalar temperature, the displacement of elastic structure and the elastic structure scalar temperature respectively, $\varepsilon>0$ denotes the viscosity coefficient for the fluid, $\mu$ is the Lam\'{e} constant for the elastic structure, $k_{1}$ and $k_{2}$ denote the thermal diffusivity, $e$ denotes the unit vector that depends only on the coordinate system, $T>0 $ is terminal time, and $f_{1}, f_{2}, f_{3}, f_{4}$ are given external force functions. The fluid domain $\Omega_{f}$ and the elastic structure domain $\Omega_{s}$ are in a fixed connected bounded domain $\Omega\subset \mathbb{R}^{\mathrm{d}}$($\mathrm{d}=2,3$) with $\Omega=\overline{\Omega_{f}}\bigcup\Omega_{s}$, $\Omega_{f}\cap\Omega_{s}=\emptyset$, the smooth interface boundary $\Gamma=\partial\Omega_{f}\cap\partial\Omega_{s}$ between the fluid and the elastic structure and the smooth boundary $\Gamma_{out}=\partial\Omega_{s}\backslash\Gamma$ satisfying $\Gamma_{out}\bigcap\Gamma=\emptyset$. The vectors $n$ and $N $ denote  the unit outer normal ones of bounded domain $\Omega_{f}$ and $\Omega_{s}$ respectively.
The model \eqref{1.1} is called a nonlinear incompressible Boussinesq-structure interaction system which describes the motion of the incompressible Boussinesq fluid surrounded by the elastic structure with the heat transfer exchange. $\eqref{1.1}_{1}$-$\eqref{1.1}_{3}$ represent incompressible Boussinesq equations, $\eqref{1.1}_{4}$-$\eqref{1.1}_{5}$ represent the linear elasticity equation with thermal diffusion. We prescribe physically coupling conditions at the interface for velocity and temperature by $\eqref{1.1}_{6}$-$\eqref{1.1}_{9}$, where 
$\eqref{1.1}_{6}$ and $\eqref{1.1}_{8}$ represent the continuity of the velocity and heat exchange at the interface $\Gamma$, $\eqref{1.1}_{7}$ and $\eqref{1.1}_{9}$ represent the stress matching conditions. $\eqref{1.1}_{10}$ is the outer boundary condition for the elastic structure on $\Gamma_{out}$. $\eqref{1.1}_{11}$-$\eqref{1.1}_{12}$ denote initial conditions for the fluid and the elastic structure respectively. Also, when $\rho=\theta=0$, the system \eqref{1.1} comes back to the classical fluid-structure interaction(FSI) model with application in biology \cite{refC}, which was first studied by J. L. Lions in \cite{refLions}, where J. L. Lions proved the existence of global weak solutions in two and three-dimensional cases as well as the uniqueness in two-dimensional case. Later, some results on the regularity of the weak solution to the system \eqref{1.1} in the case of $\rho=\theta=0$ were established by V.barel et al. \cite{refBarbu1} for the local existence of the strong solution under the assumption that $(u_{0}, w_{0}, w_{1})\in H^{2}(\Omega_{f})\times H^{2}(\Omega_{s})\times H^{1}(\Omega_{s})$, I. Kukavican et al. \cite{refKukavica1} for the local existence of the strong solution under the assumption of weaker regularity on initial data than that of V. barel et al. \cite{refBarbu1} and Shen et al \cite{refShen1} the global existence of the smooth solution in two-dimension.

There are a large number of literatures on fluid-structure interaction(FSI) problems concerning a fixed interface and a moving interface based on different physical backgrounds with many applications, for example in geology \cite{refLions}, biology \cite{refC}, biomedicine \cite{refBodnar} and aeroelasticity \cite{refDowell}.

For the FSI problem with a fixed interface for the motion of fluid inside the elastic solid which is controlled completely by infinitesimal displacements, that's only the high frequency and small displacement oscillation.
When the fluid scalar temperature and the elastic structure scalar temperature in the Boussinesq-structure interaction model \eqref{1.1} simultaneously disappear, the system \eqref{1.1} transforms into the classical FSI model, which was first introduced by J.L. Lions who in \cite{refLions} proved the existence of global weak solutions in two and three-dimensional cases as well as the uniqueness in two-dimensional space. Later, Q. Du et al. \cite{refDu} studied the well-posedness of the FSI system between a stokes fluid and a linear elastic structure with additional regularity and compatibility assumptions. For an incompressible fluid surrounded by elastic structure, V.barel et al. \cite{refBarbu1} obtained the existence and uniqueness of local smooth solution under the regular initial datum $(u_{0}, w_{0}, w_{1})\in H^{2}(\Omega_{f})\times H^{2}(\Omega_{s})\times H^{1}(\Omega_{s})$. Recently, I. Kukavican et al. \cite{refKukavica1} established the existence of local in-time smooth solutions by using an important Hidden Regularity Theorem to wave equation. The global wellposedness of the classical FSI model in two-dimensional spaces has been proved in \cite{refShen1}. For an elastic structure immersed in an incompressible fluid, Shen et al. in \cite{refShen} proved the existence of local-in-time strong solutions to the FSI model where the elastic equation with St.Venant-Kirchoff elastic material. Similarly, other related FSI models about steady interface can be founded in \cite{refBoffi,refBarbu1,refBarbu2,refCrosetto,refKukavica,refKuttler,refLx} and the references therein.\\

However, the fluid-structure interaction system with a moving interface usually be studied in Eulerian-Lagrangian coordinates. da Veiga et al. \cite{refHB} proved the existence of a locally in-time strong solution to a coupled fluid-structure evolution problem in 2D dimensional cases with periodic boundary conditions. J. Lequeurre in \cite{refLequeurre1} studied the structure modeled by a clamped viscoelastic beam and established the existence of strong solutions. D. Coutand and S. Shkoller in \cite{refCoutand1,refCoutand2} talked about the motion of an elastic solid inside an incompressible viscous fluid and proved the existence and uniqueness of the local strong solutions by applying the context of Lagrangian coordinates framework and the Tychonoff fixed-point Theorem.  Recently, there exists abundant literature about the moving interface problem that an elastic structure inside the incompressible or compressible viscous fluid, we refer to \cite{refBreit,refBae,refChambolle,refGrandmont,refLengeler1, refLasiecka} and further references.

The main purpose of this paper is to establish the global existence of weak solutions to the system \eqref{1.1} in two and three dimensions, the uniqueness of the weak solutions to the system \eqref{1.1} in two-dimension, and the existence of smooth solutions to the system \eqref{1.1} in two-dimension with smooth large initial data satisfying suitable compatibility conditions.

For simplicity, we give some notations
\begin{equation}
\begin{split}\nonumber
&(u,v)_{\Omega}=\int_{\Omega}u\cdot v dx=\int_{\Omega_{f}}u\cdot vdx+\int_{\Omega_{s}}u\cdot vdx,\\
&a_{\Omega_{i}}(u,v)=\int_{\Omega_{i}}\nabla u\cdot \nabla vdx,\quad b_{\Omega_{i}}(u,v,w)=\int_{\Omega_{i}}u\cdot \nabla v\cdot wdx,\; i=f,s.
\end{split}
\end{equation}
We introduce the following function spaces: $H^{1}(\Omega)$ ($H^{1}(\Omega_f)$ or $H^{1}(\Omega_s)$) is the standard Sobolev space, $V'$ denotes the dual space of $V$, $\mathcal{D'}$ denotes the dual space of $C_{0}^{\infty}$, and
\begin{equation}
\begin{split}\nonumber
&V_{f}=\{v(x)\in H^{1}(\Omega_{f})|\text{div}v=0, x\in\Omega_{f}\},\\
&V_{s}=\{v(x)\in H^{1}(\Omega_{s})|v|_{\Gamma_{out}}=0\},\\
&V=\{v(x)\in H_{0}^{1}(\Omega)|\Omega=\overline{\Omega}_{f}\cup \Omega_{s}, \text{div} v=0 \text{ for } x\in\Omega_{f}\},\\
&\Theta=\{\theta(x)\in H^{1}(\Omega)|\Omega=\overline{\Omega}_{f}\cup \Omega_{s}, \frac{\partial \theta}{\partial N}=0 \text{ for } x\in\Gamma_{out}\},
\end{split}
\end{equation}
and
\begin{equation}
\begin{split}\nonumber
&X_{T}=\{u\in L^{\infty}(0,T;L^{2}(\Omega))|u\in L^{2}(0,T;V_{f}),\int_{0}^{t}u d\sigma\in L^{\infty}(0,T; V_{s})\},\\
&Y_{T}=\{\rho\in L^{\infty}(0,T;L^{2}(\Omega))|\rho\in L^{2}(0,T;H^{1}(\Omega)),\frac{\partial\rho}{\partial N}=0,x\in\Gamma_{out}\}.
\end{split}
\end{equation}
Also, denote $\Phi=w_{t}$, then we have
$w(x,t)=w_{0}+\int_{0}^{t}\Phi d\sigma$ and the equation $\eqref{1.1}_{4}$ can be represented as
\begin{equation}\label{1.2}
\Phi_{t}-\mu\Delta\int_{0}^{t}\Phi d\sigma=\mu\Delta w_{0}+f_{3}.
\end{equation}
Here and in the following, $\int_0^t\Phi d\sigma=\int_0^t\Phi(\cdot,\sigma) d\sigma$.

In order to study the solution of the Boussinesq-fluid-structure interaction system \eqref{1.1}-\eqref{1.2}, we need to give a new representation way.\\
Let
\begin{equation}\nonumber
v(x, t)=\left\{\begin{aligned}
u(x, t),x\in \Omega_f\\
\Phi(x, t), x\in \Omega_s
\end{aligned}\right.,\quad
d(x, t)=\left\{\begin{aligned}
\rho(x, t), x\in \Omega_f\\
\theta(x, t), x\in \Omega_s
\end{aligned},\right.
\end{equation}
and
\begin{equation}\nonumber
f(x, t)=\left\{\begin{aligned}
f_{1}(x, t),x\in \Omega_f\\
f_{3}(x, t), x\in \Omega_s
\end{aligned}\right.,\quad
g(x, t)=\left\{\begin{aligned}
f_{2}(x, t), x\in \Omega_f\\
f_{4}(x, t), x\in \Omega_s
\end{aligned}.\right.
\end{equation}
We know that $v(x, t)$ and $d(x, t)$ are the continuous functions on $\Gamma$ because of $u(x, t)=\Phi(x, t)$, $\rho(x, t)=\theta(x, t)$ on $\Gamma$.
Then the system \eqref{1.1}-\eqref{1.2} can be rewritten into following system
\begin{align}\label{1.3}
\left \{
\begin{array}{lr}
   v_{t}+ v\cdot\nabla v+\nabla p=\varepsilon\Delta v+d e+f,&(x,t)\in\Omega_{f}\times(0,T),\\
  \text{div} v=0,&(x,t)\in\Omega_{f}\times(0,T),\\
  d_{t}+v\cdot \nabla d=k_{1}\Delta d+g,&(x,t)\in\Omega_{f}\times(0,T),\\
  v_{t}-\mu\Delta\int_{0}^{t}vd\sigma=\mu\Delta w_{0}+f,&(x,t)\in\Omega_{s}\times(0,T),\\
  d_{t}=k_{2}\Delta d+g,&(x,t)\in\Omega_{s}\times(0,T),\\
  \varepsilon\frac{\partial v}{\partial n}-pn-\frac{1}{2}(v\cdot n)v=\mu\frac{\partial}{\partial n}(\int_{0}^{t}vd\sigma+w_{0}),&(x,t)\in\Gamma\times(0,T),\\
  k_{1}\frac{\partial d}{\partial n}-\frac{1}{2} (v\cdot n)d=k_{2}\frac{\partial d}{\partial n},&(x,t)\in\Gamma\times(0,T),\\
  w_{0}+\int_{0}^{t}v d\sigma=0,\,\frac{\partial d}{\partial N}=0, &(x,t)\in\Gamma_{out}\times(0,T),\\
\end{array}
\right.
\end{align}
and initial conditional
\begin{equation}\nonumber
v(0)=v_{0}(x)=\left\{\begin{aligned}
u(0)=u_{0},x\in \Omega_f\\
w_{t}(0)=w_{1}, x\in \Omega_s
\end{aligned}\right.,\quad
d(0)=d_{0}(x)=\left\{\begin{aligned}
\rho(0)=\rho_{0}, x\in \Omega_f\\
\theta(0)=\theta_{0}, x\in \Omega_s
\end{aligned}.\right.
\end{equation}

Before stating our main results, let us first introduce the definitions of weak solutions to the system \eqref{1.3} and recall the definition of difference quotient and some important lemmas.

\begin{defn}\label{defn1.1.} The function $\{v,d\}$ defined on $\Omega\times(0,T)$ is called a weak solution to the system \eqref{1.3} if for $v\in X_{T}$,\, $d\in Y_{T}$ and
\begin{eqnarray}\label{1.4}
\begin{split}
&(v_{t}, \varphi)_{\Omega}+b_{\Omega_{f}}(v,v,\varphi)+\varepsilon a_{\Omega_{f}}(v,\varphi)+ \mu a_{\Omega_{s}}(\int_{0}^{t}v d\sigma,\varphi)\\
-&\frac{1}{2}\int_{\Gamma}(v\cdot n)v\cdot \varphi d\Gamma
=(d e,\varphi)_{\Omega_{f}}+(f,\varphi)_{\Omega}-\mu a_{\Omega_{s}}(w_{0},\varphi)\\
\end{split}
\end{eqnarray}
and
\begin{eqnarray}\label{1.5}
\begin{split}
&(d_{t},\phi)_{\Omega}+k_{1} a_{\Omega_{f}}(d,\phi)+k_{2}a_{\Omega_{s}}(d,\phi)+b_{\Omega_{f}}(v,d,\phi)\\
-&\frac{1}{2}\int_{\Gamma}(v\cdot n)d \cdot \phi d\Gamma
=(g,\phi)_{\Omega}\\
\end{split}
\end{eqnarray}
and
\begin{equation}\nonumber
v_{0}(x)=\left\{\begin{aligned}
u_{0},x\in \Omega_f\\
w_{1}, x\in \Omega_s
\end{aligned}\right.,\quad
d_{0}(x)=\left\{\begin{aligned}
\rho_{0}, x\in \Omega_f\\
\theta_{0}, x\in \Omega_s
\end{aligned}.\right.
\end{equation}
where test functions $\varphi\in V$, $\phi \in \Theta$ and initial date $v_{0}(x)$ and $d_{0}(x)$ are still the continuous functions on $\Gamma$ . Here $\mathrm{d}=2, 3$, $d\Gamma$ denotes the arc differential or area differential respectively.
\end{defn}
\textbf{Remark} We can deduce the system $\eqref{1.3}$ in some weak sense from the definition \ref{defn1.1.} as follows.
(1)Taking the  appropriate compact supported functions $\varphi$ satisfies $\text{div}\varphi=0$ in $\Omega_{f}$ and $\varphi=0$ in $\Omega_{s}$, we can get $\eqref{1.3}_{1}$ from $\eqref{1.4}$. At the same time, letting  $\varphi \in \mathcal{C}_{c}^{\infty}(\Omega)$ and $\varphi=0$ in $\Omega_{f}$, we can get $\eqref{1.3}_{4}$ from $\eqref{1.4}$. Similarly, we can easily derive $\eqref{1.3}_{3}$ and $\eqref{1.3}_{5}$ through $\eqref{1.5}$.\\
(2) Multiplying $\eqref{1.3}_{1}$ and $\eqref{1.3}_{3}$ by $\varphi$ respectively, performing integration by parts and summing the resulting equations, with the help of the continuous of $v$ on $\Gamma$, we have
\begin{equation}\label{1.6}
\begin{split}
(v_{t},\varphi)_{\Omega}&+\varepsilon a_{\Omega_{f}}(v,\varphi)+ \mu a_{\Omega_{s}}(\int_{0}^{t}v d\sigma,\varphi)+b_{\Omega_{f}}(v,v,\varphi)\\
&+\int_{\Gamma}(-\varepsilon \frac{\partial v}{\partial n}+pn+\mu\frac{\partial }{\partial n}(\int_{0}^{t}v d\sigma+w_{0}))\cdot \varphi d\Gamma\\
=&(d e,\varphi)_{\Omega_{f}}-\mu a_{\Omega_{s}}(w_{0},\varphi)+(f,\varphi)_{\Omega}.\\
\end{split}
\end{equation}
Similarly, we have 
\begin{equation}\label{1.7}
\begin{split}
(d_{t}, \phi)_{\Omega}&+k_{1} a_{\Omega_{f}}(d,\phi)+ k_{2}a_{\Omega_{s}}(d,\phi)+b_{\Omega_{f}}(v,d,\phi)\\
&+\int_{\Gamma}(k_{2}\frac{\partial d}{\partial n}-k_{1}\frac{\partial d}{\partial n})\phi d\Gamma=(g,\phi)_{\Omega}.\\
\end{split}
\end{equation}
Combining $\eqref{1.4}$ and $\eqref{1.6}$ and then combining $\eqref{1.5}$ and $\eqref{1.7}$, we get
\begin{equation}\label{1.8}
\int_{\Gamma}[-\varepsilon \frac{\partial v}{\partial n}+pn+\mu\frac{\partial }{\partial n}(\int_{0}^{t}v d\sigma+w_{0})+\frac{1}{2}(v\cdot n)v]\cdot \varphi d\Gamma=0
\end{equation}
and 
\begin{equation}\label{1.9}
\int_{\Gamma}[k_{2}\frac{\partial d}{\partial n}-k_{1}\frac{\partial d}{\partial n}+\frac{1}{2} (v\cdot n)d]\phi d\Gamma=0.
\end{equation}
Because of $\nabla\cdot \varphi=0$ in $\Omega_{f}$, we have $\int_{\Gamma}\varphi\cdot nd\Gamma=0$, which, together with \eqref{1.8}, yields that there exists a constant $\lambda$ satisfying
\begin{equation}\label{1.10}
-\varepsilon \frac{\partial v}{\partial n}+pn+\mu\frac{\partial}{\partial n}(\int_{0}^{t}v d\sigma+w_{0})+\frac{1}{2}(v\cdot n)v=\lambda n.
\end{equation}
Letting $p-\lambda$ instead of $p$ in $\eqref{1.10}$, we get $\eqref{1.3}_{6}$. Also, it is easy to get $\eqref{1.3}_{7}$ from $\eqref{1.9}$. Then, we can get the rest equations or the boundary conditions.

\begin{defn}\label{defn1.2.}(\cite{refCoutand1})
The first-order difference quotient $D_{h}u(x_{1},x_{2},t)$ and the second-order difference quotient $D_{-h}D_{h}u(x_{1},x_{2},t)$ are defined respectively by
\begin{align}\nonumber
\begin{split}
&D_{h}u(x_{1},x_{2},t)=\frac{u(x_{1}+h,x_{2},t)-u(x_{1},x_{2},t)}{h},\\
&D_{-h}D_{h}u(x_{1},x_{2},t)=\frac{u(x_{1}+h,x_{2},t)-2u(x_{1},x_{2},t)+u(x_{1}-h,x_{2},t)}{h^{2}},\\
\end{split}
\end{align}
where $h\in \mathbb{R}$. Hence, the following properties hold
\begin{align}\label{1.11}
\begin{split}
&(u,D_{-h}D_{h}u)_{\Omega}=(D_{h}u,D_{h}u)_{\Omega},\\
&D_{h}(uv)=hD_{h}uD_{h}v+uD_{h}v+vD_{h}u.\\
\end{split}
\end{align}
\end{defn}
\textbf{Lemma 1.1}(\cite{refLions1}) Let $\Omega_{f} \subset\mathbb{R}^{\mathrm{d}}$ be a bounded domain with the smooth boundary $\Gamma=\partial{\Omega_{f}}$. If $v\in H^{1-\epsilon}(\Omega_{f})$ and $u\in H^{1}(\Omega_{f})$, then
\begin{equation}\label{1.12}
\begin{split}
&\|v\|_{L^{3}(\Omega_{f})}\leq \|v\|_{H^{1/2}(\Omega_{f})}\leq \|v\|_{H^{1-\epsilon}(\Omega_{f})},\, 0<\epsilon<\frac{1}{2},\, \mathrm{d}=2, 3,\\
&\|u|_{\Gamma}\|_{L^{3}(\Gamma)}\leq C\|u\|^{\tfrac{\mathrm{d}+2}{6}}_{H^{1}(\Omega_{f})}\|u\|^{1-\tfrac{\mathrm{d}+2}{6}}_{L^{2}(\Omega_{f})}, \mathrm{d}=2,3,\\
&\|u|_{\Gamma} \|_{L^{4}(\Gamma)}\leq C \|u\|_{H^{1/2}(\Gamma)}\leq \|u\|_{H^{1}(\Omega_{f})},\, \mathrm{d}=2, 3.\\
\end{split}
\end{equation}
\textbf{Lemma 1.2}(\cite{refLadyzhenskaya}) Let $\Omega_{f} \subset\mathbb{R}^{\mathrm{d}}$ be a bounded domain with the smooth boundary $\Gamma=\partial{\Omega_{f}}$. If $u\in H^{1}(\Omega_{f})$, then
\begin{equation}\label{1.13}
\|u\|_{L^{4}(\Omega_{f})}\leq C\|u\|_{H^{1}(\Omega_{f})}^{\tfrac{\mathrm{d}}{4}}\|u\|_{L^{2}(\Omega_{f})}^{1-\tfrac{\mathrm{d}}{4}}, \mathrm{d}=2,3.
\end{equation}
\textbf{Lemma 1.3}(\cite{refBoyer}) Let $1\le p_0, p_1\le +\infty$. The spaces $B_{0}$, $B$ and $B_{1}$ are Banach spaces and satisfy the following conditions:
(i)$\{u_{i}\}_{i=1}^{\infty}$ is bounded in $L^{p_{0}}(0,T;B_{0})$,
(ii)$\{u_{i,t}\}_{i=1}^{\infty}$ is bounded in $L^{p_{1}}(0,T;B_{1})$, and
(iii)$B_{0}\circlearrowleft \circlearrowleft B \circlearrowleft B_{1}$.
Then
there exists a strongly converging subsequence of $\{u_{i}\}_{i=1}^{\infty}$ in $L^{p_{0}}(0,T;B)$ provided $p_{0}<\infty$,
and there exists a strongly converging subsequence of  $\{u_{i}\}_{i=1}^{\infty}$ in $C(0,T;B)$ provided $p_{0}=\infty$ and $p_{1}>1$.\\
\\
Now, we state our main results as follows.
\begin{thm}\label{thm 1.1.} Let $\mathrm{d}=2 \;{\rm or }\; 3$. Assume that the domain $\Omega_{f}$, $\Omega_{s}$ are smooth and bounded with the smooth interface boundary $\Gamma=\partial\Omega_{f}\cap\partial\Omega_{s}$. And for the initial date satisfy $(v_{0},d_{0})\in L^{2}(\Omega)$, $w_{0}\in V_{s}$, $\text{div} v_{0}=0$ in $\Omega_{f}$, $(v_0(x), d_0(x))\in C$  on $\Gamma$ and $\frac{\partial d_0}{\partial N}=0$  on $\Gamma_{out}$. Also, assume that $(f_{1},f_{2})\in L^{2}(0, \infty; L^{2}(\Omega_{f}))$ and $(f_{3},f_{4})\in L^{2}(0, \infty; L^{2}(\Omega_{s}))$. Then, there exists a global weak solution of \eqref{1.3} in the sense of the definition \ref{defn1.1.}, defined in the internal $[0,T]$. Also, it satisfies, for some $1<p<2$,
\begin{equation}\nonumber
v_{t} \in L^{p}(0,T; V'),\; d_{t}\in L^{p}(0,T;\Theta').
\end{equation}
Moreover, when $\mathrm{d}=2$, the global weak solution of \eqref{1.3} is unique.
\end{thm}

Furthermore, if the initial data satisfy more assumptions on the compatibility conditions in two-dimensional case, we have the following results on the strong and smooth solutions.

\begin{thm}\label{thm 1.2.} Let $\mathrm{d}=2$, $L$ is a fixed constant, $\Omega_{f}=\mathcal{T}\times(-\tfrac{L}{2}, \tfrac{L}{2})$, $\Omega_{s}=\mathcal{T}\times((-L,-\tfrac{L}{2})\cup(\tfrac{L}{2},L))$ with the torus $\mathcal{T}=\tfrac{\mathbb{R}}{2\pi\mathcal{Z}}$. In this case, $\Gamma=\{(x,y)\in\mathbb{R}^{2}|x\in\mathcal{T},\;y=-\tfrac{L}{2}\;and\; y=\tfrac{L}{2}\}$ and $\Gamma_{out}=\{(x,y)|x\in\mathcal{T},\;y=-L\; and\;y=L\}$, and $n=N=\{0, \pm 1\}$. The operator $\Pi$ denotes the projection on the tangent space of the boundary $\Gamma$ i.e, $\Pi=I\mathrm{d}-n\otimes n$. Assume that the initial data $(\rho_{0}, u_{0},w_{0},w_{1},\theta_{0})$ satisfy all assumptions given by Theorem \ref{thm 1.1.} and, furthermore, satisfy
\begin{equation}\label{1.14}
\begin{cases}
u_{0}\in V_{f}\cap H^{2}(\Omega_{f}),\rho_{0} \in H^{2}(\Omega_{f}),\\
w_{0}\in V_{s}\cap H^{2}(\Omega_{s}), w_{1}\in V_{s}, \theta_{0}\in H^{2}(\Omega_{s})\\
(f_{1}, f_{2})\in  H^{1}(0, \infty; L^{2}(\Omega_{f}))\cap L^{2}(0, \infty; H^{1}(\Omega_{f})),\\
(f_{3}, f_{4})\in H^{1}(0, \infty; L^{2}(\Omega_{s}))\cap L^{2}(0, \infty; H^{1}(\Omega_{s}))
\end{cases}
\end{equation}
and the zero order compatibility conditions
\begin{large}
\begin{equation}\label{1.15}
\begin{cases}
 \Pi\left(\mu\frac{\partial w_{0}}{\partial n}-\varepsilon\frac{\partial u_{0}}{\partial n}+\frac{1}{2}(u_{0}\cdot n)u_{0}\right)=0, x\in\Gamma,\\
 k_{2}\frac{\partial \theta_{0}}{\partial n}-k_{1}\frac{\partial \rho_{0}}{\partial n}+\frac{1}{2}(u_{0}\cdot n)\rho_{0}=0,x\in\Gamma.
\end{cases}
\end{equation}
\end{large}
Then the global weak solution $\{v,d\}$ of system \eqref{1.3} become the global strong solution satisfying
\begin{equation}\label{1.16}
\begin{cases}
v_{t}\in X_{T},d_{t}\in Y_{T},\\
v\in L^{\infty}(0,T; V_{f})\cap L^{2}(0,T; H^{2}(\Omega_{f})),\\
v\in L^{\infty}(0,T;V_{s}),\int_{0}^{t}v d \sigma \in L^{\infty}(0,T;H^{2}(\Omega_{s})),\\
p\in L^{2}(0,T; H^{1}(\Omega_{f})),\\
d\in L^{\infty}(0,T; H^{1}(\Omega))\cap L^{2}(0,T; H^{2}(\Omega)).
\end{cases}
\end{equation}
\end{thm}
Moreover, if the initial value $(v_{0}, d_{0})$ and the external force function $(f_{1},f_{2},f_{3},f_{4})$ are smooth and satisfy suitable high order compatibility conditions at the interface, then the global weak solution is the global smooth solution.\\

Next, we will show the global existence in two/three-dimensional and uniqueness in two-dimensional of the weak solution of \eqref{1.1} in section 2 and the global existence of strong solution of \eqref{1.1} in section 3.

\section{\bf Proof of Theorem 1.1}
The proof of Theorem 1.1 contains four steps. Step 1 is to construct solutions of finite-dimensional approximation to the system \eqref{1.3} by the well-known Galerkin method. Here we should point out that the existence of the basis of the solution's space here can be guaranteed by the fact that the separable Hilbert space has one complete countable orthogonal basis, but we have not their exact expression form, and, hence, some technical difficulties are required to deal with. Step 2 is to get the uniform estimates of the finite-dimensional approximation solutions by the energy method. Step 3 is to achieve a weak solution of the system \eqref{1.3} by the compactness method, which is the main difficulty to solve nonlinear terms. Step 4 is to prove that the global weak solution is unique in two dimensions.\\
\\
\textbf{Step 1: Galerkin approximations}\\

Through a finite-dimensional approximation, we can construct a weak solution.
Let $\{\varphi_{j}\}_{j=1}^{\infty}$ and $\{\phi_{j}\}_{j=1}^{\infty}$ 
be an orthogonal basis of $V$ and $\Theta$ satisfying $\|\varphi_j\|_{L^2(\Omega)}=\|\phi_j\|_{L^2(\Omega)}=1$ respectively, where the existence of the orthogonal basis for the space $V$ and the space $\Theta$ is because they are all separable Hilbert spaces.

Fixing a positive integer $m$, we will find a class of approximating solution $v_{m}$: $[0,T]\rightarrow V$ and  $d_{m}$: $[0,T]\rightarrow \Theta$ of the form
\begin{equation}\label{2.1}
v_{m}(x,t)=\sum_{j=1}^{m}a_{j}^{m}(t)\varphi_{j}(x),\quad d_{m}(x,t)=\sum_{j=1}^{m}b_{j}^{m}(t)\phi_{j}(x)
\end{equation}
and select the vector $a_{j}^{m}(t)$ and $b_{j}^{m}(t)$ for $ 0\leq t\leq T$$(j\leq m,m=1,2\cdots)$ such that
they satisfies the following ordinary differential equations(ODEs)
\begin{equation}\label{2.2}
\begin{split}
&(v_{m,t},\varphi_{j})_{\Omega}+\varepsilon a_{\Omega_{f}}(v_{m},\varphi_{j})+\mu a_{\Omega_{s}}(\int_{0}^{t}v_{m} d\sigma,\varphi_{j})
+b_{\Omega_{f}}(v_{m},v_{m},\varphi_{j})\\
-&\frac{1}{2}\int_{\Gamma}(v_{m}\cdot n)v_{m}\cdot \varphi_{j}d\Gamma=(f,\varphi_{j})_{\Omega}+(d_{m}e,\varphi_{j})_{\Omega_{f}}-\mu a_{\Omega_{s}}(w_{0m},\varphi_{j}),\\
&(d_{m,t},\phi_{j})_{\Omega}+k_{1} a_{\Omega_{f}}(d_{m},\phi_{j})+k_{2}a_{\Omega_{s}}(d_{m},\phi_{j})
+b_{\Omega_{f}}(v_{m},d_{m},\phi_{j})\\
&-\frac{1}{2}\int_{\Gamma} (v_{m}\cdot n)d_{m}\cdot \phi_{j}d\Gamma=(g,\phi_{j})_{\Omega},\\
&j=1,2,\cdots,m,\;t\in[0,T]
\end{split}
\end{equation}
and the following initial conditions
\begin{equation}\label{2.3}
\begin{split}
&v_{m}(x,t=0)=\sum_{j=1}^{m}\alpha_{j}^{m}\varphi_{j}(x)\rightarrow v_{0}(x),\text{strongly in}\,V, \\
&w_{0m}(x,t=0)=\sum_{j=1}^{m}\vartheta_{j}^{m}\varphi_{j}(x)\rightarrow w_{0}(x),\text{strongly in} \,V_{s}, \\
&d_{m}(x,t=0)=\sum_{j=1}^{m}\beta_{j}^{m}\phi_{j}(x)\rightarrow d_{0}(x),\text{strongly in} \,\Theta,
\end{split}
\end{equation}  \\
as $m\rightarrow \infty$.

According to the existence and uniqueness theory for ODEs, there exists a unique continuous function $\{a_{j}^{m}(t), b_{j}^{m}(t)\}$, $j=1,\cdots,m$ satisfying $\eqref{2.2}-\eqref{2.3}$ for $0\leq t\leq t_{m}:=t(m)$. Next, we can prove $t(m)=T$ by energy estimates.\\
\\
\textbf{Step 2: Energy estimates}\\
\textbf{Lemma 2.1.} Under the assumptions of Theorem \ref{thm 1.1.}, there exists a positive constant $C(T)$ (independent of $m$) depending on the initial data $(v_{0},d_{0})$, $\Omega_{f}$, $\Omega_{s}$, $\|f\|_{L^{2}(0, \infty; L^{2}(\Omega))}$, $\|g\|_{L^{2}(0, \infty; L^{2}(\Omega))}$ and $T$, such that
\begin{equation}\label{2.4}
\begin{split}
&\|v_{m}\|_{L^{\infty}(0,T; {L^{2}(\Omega))}}^{2}
+\|d_{m}\|_{L^{\infty}(0,T; {L^{2}(\Omega))}}^{2}\\
+&\|\int_{0}^{t}v_{m}d\sigma+w_{0m}\|_{L^{\infty}(0,T; V_{s})}^{2}
+\| v_{m}\|_{L^{2}(0,T;V_{f})}^{2}+ \|d_{m}\|_{L^{2}(0,T;H^{1}(\Omega))}^{2}
\leq C(T).
\end{split}
\end{equation}
\textbf{Proof.} Taking the inner product of $\eqref{2.2}_{1}$, $\eqref{2.2}_{2}$ with $a_{j}^{m}(t)$, $b_{j}^{m}(t)$ respectively and summing with respect to $j=1,2,\cdot\cdot\cdot,m$, we have
\begin{equation}\label{2.5}
\begin{split}
\frac{1}{2}\frac{d}{dt}&(\|v_{m}\|_{L^{2}(\Omega)}^{2}+\mu\|\nabla (\int_{0}^{t}v_{m}d\sigma+w_{0m})\|_{L^{2}(\Omega_{s})}^{2})+\varepsilon\|\nabla v_{m}\|_{L^{2}(\Omega_{f})}^{2}\\
+&b_{\Omega_{f}}(v_{m},v_{m},v_{m})
-\frac{1}{2}\int_{\Gamma}(v_{m}\cdot n)|v_{m}|^{2}d\Gamma=(f,v_{m})_{\Omega}+(d_{m}e,v_{m})_{\Omega_{f}},\\
\end{split}
\end{equation}
and
\begin{equation}\label{2.6}
\begin{split}
&\frac{1}{2}\frac{d}{dt}(\|d_{m}\|_{L^{2}(\Omega)}^{2})+k_{1} \|\nabla d_{m}\|_{L^{2}(\Omega_{f})}^{2}+k_{2} \|\nabla d_{m}\|_{L^{2}(\Omega_{s})}^{2}\\
&+b_{\Omega_{f}}(v_{m},d_{m},d_{m})-\frac{1}{2}\int_{\Gamma}( v_{m}\cdot n)| d_{m}|^{2}d\Gamma=(g,d_{m})_{\Omega},\\
\end{split}
\end{equation}
where
\begin{equation}\nonumber
\begin{split}
&(f,v_{m})_{\Omega}+(d_{m}e,v_{m})_{\Omega_{f}}\\
\leq& \|f\|_{L^{2}(\Omega)}\|v_{m}\|_{L^{2}(\Omega)}+\|d_m\|_{L^{2}(\Omega)}\|v_{m}\|_{L^{2}(\Omega)}\\
\leq& \|v_{m}\|^{2}_{L^{2}(\Omega)}+\|f\|_{L^{2}(\Omega)}^{2}+\|d_{m}\|^{2}_{L^{2}(\Omega)}.
\end{split}
\end{equation}
Thanks to $\text{div}v_{m}=0$ in $\Omega_{f}$, we know that
\begin{equation}\label{2.7}
\begin{split}
&b_{\Omega_{f}}(v_{m},v_{m},v_{m})=\frac{1}{2}\int_{\Omega_{f}}v_{m}\cdot \nabla |v_{m}|^{2}dx=\frac{1}{2}\int_{\Gamma}|v_{m}|^{2}n\cdot v_{m}d\Gamma,\\
&b_{\Omega_{f}}(v_{m},d_{m},d_{m})=\frac{1}{2}\int_{\Omega_{f}}v_{m}\cdot \nabla |d_{m}|^{2}dx=\frac{1}{2}\int_{\Gamma}|d_{m}|^{2}n\cdot v_{m}d\Gamma.
\end{split}
\end{equation}
Combining the above estimates $\eqref{2.5}-\eqref{2.7}$, implies
\begin{equation}\label{2.8}
\begin{split}
&\frac{1}{2}\frac{d}{dt}(\|v_{m}\|_{L^{2}(\Omega)}^{2}+\|d_{m}\|_{L^{2}(\Omega)}^{2}
+\mu\|\nabla(\int_{0}^{t}v_{m}d\sigma+w_{0m})\|_{L^{2}(\Omega_{s})}^{2})\\
&+\varepsilon \|\nabla v_{m}\|_{L^{2}(\Omega_{f})}^{2}+k_{1} \|\nabla d_{m}\|_{L^{2}(\Omega_{f})}^{2}+k_{2} \|\nabla d_{m}\|_{L^{2}(\Omega_{s})}^{2}\\
\leq &\|v_{m}\|_{L^{2}(\Omega)}^{2}+\|d_m\|_{L^{2}(\Omega)}^{2}+\|f\|_{L^{2}(\Omega)}^{2}+\|g\|_{L^{2}(\Omega)}^{2},\quad t\in[0,T].
\end{split}
\end{equation}\\
Applying Gronwall's inequality to \eqref{2.8} and using the assumptions in Theorem \ref{thm 1.1.}, we have
 \begin{equation*}
\begin{split}
&\|v_{m}(t)\|_{L^{2}(\Omega)}^{2}+\|d_{m}(t)\|_{L^{2}(\Omega)}^{2}
+\mu\|\nabla (\int_{0}^{t}v_{m}d\sigma+w_{0m})\|_{L^{2}(\Omega_{s})}^{2}\\
&+\int_{0}^{t}(\varepsilon \|\nabla v_{m}\|_{L^{2}(\Omega_{f})}^{2}+k_{1} \|\nabla d_{m}\|_{L^{2}(\Omega_{f})}^{2}+k_{2} \|\nabla d_{m}\|_{L^{2}(\Omega_{s})}^{2})dt\\
\leq &e^{Ct}(\|v_{0}\|_{L^{2}(\Omega)}^{2}+\|d_{0}\|_{L^{2}(\Omega)}^{2}+\mu\|\nabla w_{0}\|_{L^{2}(\Omega_{s})}^{2})\\
\leq &C(T)<\infty,\quad t\in[0,T],
\end{split}
\end{equation*}\\ which gives the desired estimate \eqref{2.4}.

By using \eqref{2.4}, we know that there exists a positive constant $C(T)$ such that $|a_{j}^{m}(t)|+| b_{j}^{m}(t)|\leq C(T)<\infty$ for any $m\geq 1$, $j=1,2,\cdots, m$ and $0\leq t\leq t(m)$, and hence, $t(m)=T$ by the extension theory of ODEs. Next, we will give a priori estimates of $\{v_{m,t}, d_{m,t}\}$.\\
\\
\textbf{Lemma 2.2.} Under the assumptions of Theorem \ref{thm 1.1.}, the following hold uniformly on $m=1,2,\cdots$
\begin{equation}\label{2.9}
v_{m,t}\in  L^{p}(0,T; V'),\; d_{m,t}\in L^{p}(0,T;\Theta').
\end{equation}
\textbf{Proof.}
For $j=1,2\cdots,m$, by the fact that $u_{m}|_{\Gamma}\neq 0$, $\varphi_{j}|_{\Gamma}\neq 0$.
The equation $\eqref{2.2}_{1}$ can be rewritten as
\begin{eqnarray}\label{2.10}
\begin{split}
&(v_{m,t},\varphi_{j})_{\Omega}+\varepsilon a_{\Omega_{f}}(v_{m},\varphi_{j})+\mu a_{\Omega_{s}}(\int_{0}^{t}v_{m}d\sigma,\varphi_{j})
+b_{\Omega_{f}}(v_{m},v_{m},\varphi_{j})\\
-&\frac{1}{2}\int_{\Gamma}(v_{m}\cdot n)v_{m}\cdot \varphi_{j} d\Gamma=(f,\varphi_{j})_{\Omega}+(d_{m}e,\varphi_{j})_{\Omega_{f}}-\mu a_{\Omega_{s}}(w_{0},\varphi_{j}),\\
\end{split}
\end{eqnarray}
where
\begin{equation}\nonumber
b_{\Omega_{f}}(v_{m},v_{m},\varphi_{j})=-b_{\Omega_{f}}(v_{m},\varphi_{j},v_{m})+\int_{\Gamma}(v_{m}\cdot n)\cdot(v_{m}\cdot \varphi_{j})d\Gamma.
\end{equation}
Using H\"{o}lder's inequality, Young's inequality, \eqref{2.4}, $\eqref{1.12}_{2}$ and \eqref{1.13} to $\eqref{2.10}$ and  the fact that $1\leq\frac{\mathrm{d}}{2}<\frac{\mathrm{d}+2}{3}<2$ $(\mathrm{d}=2,3)$, we have, for $j=1,2\cdots,m$,
\begin{flalign}\tag{2.11}\label{2.11}
|(v_{m,t},\varphi_{j})|
=&|-\varepsilon a_{\Omega_{f}}(v_{m},\varphi_{j})+(d_{m}e,\varphi_{j})_{\Omega_{f}}
+b_{\Omega_{f}}(v_{m},\varphi_{j}, v_{m})\\ \nonumber
&-\int_{\Gamma}(v_{m}\cdot n)\cdot(v_{m}\cdot \varphi_{j})d\Gamma+\frac{1}{2}\int_{\Gamma}(v_{m}\cdot n)v_{m}\cdot \varphi_{j} d\Gamma+(f,\varphi_{j})_{\Omega}|\\ \nonumber
&+|-\mu a_{\Omega_{s}}(\int_{0}^{t}v_{m}d\sigma+w_{0m},\varphi_{j})|\\ \nonumber
\leq&C(\|\nabla v_{m}\|_{L^{2}(\Omega_{f})}+\|v_{m}\|^{2}_{L^{4}(\Omega_{f})}+\|\nabla(\int_{0}^{t}v_{m}d\sigma+w_{0m})\|_{L^{2}(\Omega_{s})})\|\nabla \varphi_{j}\|_{L^{2}(\Omega)}\\ \nonumber
&+\|d_{m}\|_{L^{2}(\Omega_{f})}\|\varphi_{j}\|_{L^{2}(\Omega_{f})}+\|v_{m}|_{\Gamma}\|^{2}_{L^{3}(\Gamma)}\|\varphi_{j}\|_{L^{3}(\Gamma)}+\|f\|_{L^{2}(\Omega)}\| \varphi_{j}\|_{L^{2}(\Omega)}\\ \nonumber
\leq&C(\|v_{m}\|_{V_{f}}+\| v_{m}\|^{\tfrac{\mathbf{d}}{2}}_{V_{f}}+\|v_{m}\|^{\tfrac{\mathbf{d}+2}{3}}_{V_{f}}+\|\int_{0}^{t}v_{m}d\sigma+w_{0m}\|_{V_{s}})\| \nabla \varphi_{j}\|_{L^{2}(\Omega)}\\ \nonumber
&+C(\|d_{m}\|_{L^{2}(\Omega_{f})}+\|f\|_{L^{2}(\Omega)})\| \varphi_{j}\|_{L^{2}(\Omega)}\\ \nonumber
\leq&C(T)(\| v_{m}\|^{\tfrac{\mathbf{d}+2}{3}}_{V_{f}}+\|\int_{0}^{t}v_{m}d\sigma+w_{0m}\|_{V_{s}}+\|d_{m}\|_{L^{2}(\Omega_{f})}+\|f\|_{L^{2}(\Omega)})\|\varphi_{j}\|_{V}, \nonumber
\end{flalign}

which implies that
\begin{equation}\label{2.12}
\|v_{m,t}\|_{V'}\leq C(T)(\| v_{m}\|^{\tfrac{\mathrm{d}+2}{3}}_{V_{f}}+\|\int_{0}^{t}v_{m}d\sigma+w_{0m}\|_{V_{s}}+\|d_{m}\|_{L^{2}(\Omega_{f})}+\|f\|_{L^{2}(\Omega)}).
\end{equation}
Combining the energy estimates \eqref{2.4} and inequality \eqref{2.12}, we know that there exists a positive constant $C(T)$ such that
\begin{equation}\label{2.13}
\|v_{m,t}\|_{L^{p}(0,T;V')}\leq C(T)<\infty,\;p=\frac{6}{\mathrm{d}+2}\in(1,2),\;\mathrm{d}=2,3.
\end{equation}\\
In $\eqref{2.2}_{2}$, we consider the equality that $b_{\Omega_{f}}(u_{m},d_{m},\phi_{j})=-b_{\Omega_{f}}(u_{m},\phi_{j},d_{m})+\int_{\Gamma}(u_{m}\cdot n)d_{m}\phi_{j}d\Gamma$ and use a similar ways as \eqref{2.11}, for $j=1,2\cdots,m$ yields to
\begin{flalign}\nonumber
|(d_{m,t},\phi_{j})_{\Omega}|=&|-k_{1} a_{\Omega_{f}}(d_{m},\phi_{j})-k_{2}a_{\Omega_{s}}(d_{m},\phi_{j})
-b_{\Omega_{f}}(v_{m},d_{m},\phi_{j})\\ \nonumber
&+\frac{1}{2}\int_{\Gamma} (v_{m}\cdot n)d_{m}\cdot \phi_{j}d\Gamma+(g,\phi_{j})_{\Omega}|\\ \nonumber
\leq &C(\|\nabla d_{m}\|_{L^{2}(\Omega)}+\|v_{m}\|_{L^{4}(\Omega_{f})}\|d_{m}\|_{L^{4}(\Omega_{f})})\|\nabla\phi_{j}\|_{L^{2}(\Omega)}\\ \nonumber
&+\|v_{m}|_{\Gamma}\|_{L^{3}(\Gamma)}\|d_{m}|_{\Gamma}\|_{L^{3}(\Gamma)}\|\phi_{j}\|_{L^{3}(\Gamma)}+\|g\|_{L^{2}(\Omega)}\| \phi_{j}\|_{L^{2}(\Omega)}\\  \nonumber
\leq &C(\| d_{m}\|_{\Theta}+\|v_{m}\|^{\tfrac{\mathbf{d}}{4}}_{V_{f}}\| d_{m}\|^{\tfrac{\mathbf{d}}{4}}_{H^{1}(\Omega_{f})})\|\nabla\phi_{j}\|_{L^{2}(\Omega)}\\ \nonumber
&+\| v_{m}\|^{\tfrac{\mathbf{d}+2}{6}}_{V_{f}}\| d_{m}\|^{\tfrac{\mathbf{d}+2}{6}}_{H^{1}(\Omega_{f})}\|\nabla\phi_{j}\|_{L^{2}(\Omega)}+\|g\|_{L^{2}(\Omega)})\| \phi_{j}\|_{L^{2}(\Omega)}\\ \nonumber
\leq&C(T)(\|v_{m}\|^{\tfrac{\mathbf{d}+2}{3}}_{V_{f}}+\|d_{m}\|^{\tfrac{\mathbf{d}+2}{3}}_{\Theta}+\|d_{m}\|_{L^{2}(\Omega_{f})}+\|g\|_{L^{2}(\Omega)})\|\phi_{j}\|_{\Theta}. \nonumber
\end{flalign}
It is easily to get that
\begin{equation}\label{2.14}
\|d_{m,t}\|_{L^{p}(0,T;\Theta')}\leq C(T)<\infty,\;p=\frac{6}{\mathrm{d}+2}\in(1,2),\;\mathrm{d}=2,3.
\end{equation}
Thus, by \eqref{2.13} and \eqref{2.14}, we know that $v_{m,t}$ is bounded in $L^{p}(0,T; V')$,  and $d_{m,t}$ is bounded in $L^{p}(0,T; \Theta')$ with $p=\frac{6}{\mathrm{d}+2}\in(1,2),\; \mathrm{d}=2,3$.\\
\\
\textbf{Step 3: Limiting processes and existence of the weak solutions}\\

1. \eqref{2.4} and \eqref{2.9} imply
\begin{equation}\label{2.15}
\begin{split}
&\{v_{m}\}_{m=1}^{\infty}\;{\rm is\; bounded \;in}\; X_{T},\\
&\{d_{m}\}_{m=1}^{\infty}\;{\rm is\; bounded \;in}\; Y_{T},\\
&\{v_{m,t}\}_{m=1}^{\infty}\;{\rm is\; bounded \;in}\; L^{p}(0,T; V'),\;1<p<2,\\
&\{d_{m,t}\}_{m=1}^{\infty}\;{\rm is\; bounded \;in}\; L^{p}(0,T;\Theta'),\;1<p<2.\\
\end{split}
\end{equation}\\
Therefore, there exist subsequences $\{v_{\mu}\}_{\mu=1}^{\infty}\subset \{v_{m}\}_{m=1}^{\infty}$, $\{d_{\mu}\}_{\mu=1}^{\infty}\subset \{d_{m}\}_{m=1}^{\infty}$ and some functions
\begin{equation}\nonumber
\begin{split}
&v\in X_{T},\,d\in Y_{T},\\
&\{v_{t},d_{t}\} \in L^{p}(0,T; V')\times L^{p}(0,T;\Theta'),
\end{split}
\end{equation}
such that
\begin{equation}\label{2.16}
\begin{split}
&v_{\mu}\rightharpoonup v \;{\rm weakly }*\;{\rm in}\; L^{\infty}(0,T;L^{2}(\Omega)),\\
&v_{\mu}\rightharpoonup v \;{\rm weakly}\;{\rm in}\; L^{2}(0,T;V_{f}),\,\int_{0}^{t}v_{\mu}d\sigma\rightharpoonup \int_{0}^{t}v d\sigma\;{\rm weakly} * \;{\rm in}\; L^{\infty}(0,T;V_{s}), \\
&d_{\mu}\rightharpoonup d \;{\rm weakly} \;{\rm in}\; L^{2}(0,T; H^{1}(\Omega)), d_{\mu}\rightharpoonup d\;{\rm weakly }*\;{\rm in}\; L^{\infty}(0,T;L^{2}(\Omega)),\\
&\{v_{\mu,t},d_{\mu,t}\} \rightharpoonup \{v_{t},d_{t}\} \;{\rm weakly}\;{\rm in}\; L^{p}(0,T; V')\times L^{p}(0,T;\Theta'). \\
\end{split}
\end{equation}
Because of $\Omega_{f}\subset \Omega$, we know that
\begin{equation*}
H^{1}(\Omega_{f})\circlearrowleft \circlearrowleft H^{1}(\Omega).
\end{equation*}
Thanks to Lemma 1.3 and the compactness embedding $H^{1}(\Omega)\circlearrowleft \circlearrowleft H^{1-\epsilon}(\Omega)$, the following conclusions hold for $0<\epsilon<\frac{1}{2}$,
\begin{equation}\label{2.17}
\begin{split}
&v_{\mu}\rightarrow v \;{\rm strongly}\;{\rm in}\; L^{2}(0,T;H^{1-\epsilon}(\Omega_{f})),\\
&d_{\mu}\rightarrow d \;{\rm strongly}\;{\rm in}\; L^{2}(0,T;H^{1-\epsilon}(\Omega)).\\
\end{split}
\end{equation}

2. Both $\Omega_{f}$ and $\Omega_{s}$ are bounded and smooth domains, $v \rightharpoonup v|_{\Gamma}$ is a continuous mapping from $H^{1-\epsilon}(\Omega)$ to $L^{2}(\Gamma)$, we have\\
\begin{eqnarray}\label{2.18}
\begin{split}
&v_{\mu}|_{\Gamma}\rightarrow v|_{\Gamma}\;{\rm strongly} \; {\rm in}\; L^{2}(0,T; L^{2}(\Gamma)),\\
&d_{\mu}|_{\Gamma}\rightarrow d|_{\Gamma}\;{\rm strongly} \; {\rm in}\; L^{2}(0,T; L^{2}(\Gamma)).
\end{split}
\end{eqnarray}
In order to take the limit as $m\rightarrow \infty$ in ODEs \eqref{2.2}, let us give the convergence of the nonlinear terms as follows.\\
Firstly, let us  prove that
\begin{equation}\label{2.19}
\int_{\Gamma}(v_{\mu}\cdot n)(d_{\mu}\cdot \phi_{j})d\Gamma \rightarrow \int_{\Gamma}(v\cdot n)(d \cdot \phi_{j})d\Gamma\;\;{\rm strongly}\;{\rm in}\; L^{1}(0,T).
\end{equation}
In fact, by using H\"{o}lder's inequality, $\eqref{1.12}_{3}$, $\eqref{2.4}$ and $\eqref{2.18}$, we conclude that
\begin{equation}\label{2.20}
\begin{split}
&\int_{\Gamma}(v_{\mu}\cdot n)(d_{\mu} \cdot \phi_{j})d\Gamma-\int_{\Gamma}(v\cdot n)(d \cdot \phi_{j})d\Gamma\\
=&\int_{\Gamma}(v_{\mu}\cdot n)[(d_{\mu}-d) \phi_{j}]d\Gamma+\int_{\Gamma}[(v_{\mu}-v)\cdot n](d\cdot \phi_{j})d\Gamma\\
\leq&\|d_{\mu}-d\|_{L^{2}(\Gamma)}\|v_{\mu}\|_{L^{4}(\Gamma)}\| \phi_{j}\|_{L^{4}(\Gamma)}+\|d\|_{L^{4}(\Gamma)}\|v_{\mu}-v\|_{L^{2}(\Gamma)}\| \phi_{j}\|_{L^{4}(\Gamma)}\\
\leq&C\|d_{\mu}-d\|_{L^{2}(\Gamma)}\|v_{\mu}\|_{H^{1/2}(\Gamma)}\| \phi_{j}\|_{H^{1/2}(\Gamma)}+C\|d\|_{H^{1/2}(\Gamma)}\|v_{\mu}-v\|_{L^{2}(\Gamma)}\| \phi_{j}\|_{H^{1/2}(\Gamma)}\\
\leq&C\|d_{\mu}-d\|_{L^{2}(\Gamma)}\|v_{\mu}\|_{V_{f}}\| \phi_{j}\|_{H^{1}(\Omega_{f})}+C\|d\|_{H^{1}(\Omega_{f})}\|v_{\mu}-v\|_{L^{2}(\Gamma)}\|\phi_{j}\|_{H^{1}(\Omega_{f})}\\
& \rightarrow 0, \quad \mu\rightarrow\infty.
\end{split}
\end{equation}\\
Similarly, we have
\begin{equation}\label{2.21}
\int_{\Gamma}(v_{\mu}\cdot n)v_{\mu}\cdot \varphi_{j}d\Gamma\rightarrow\int_{\Gamma}(v\cdot n)v\cdot \varphi_{j}d\Gamma\;\;{\rm strongly}\;{\rm in}\; L^{1}(0,T).
\end{equation}\\
Secondly, let us prove
\begin{equation}\label{2.22}
\int_{\Omega_{f}}v_{\mu}\cdot\nabla d_{\mu}\cdot \phi_{j}dx\rightarrow \int_{\Omega_{f}}v\cdot\nabla d\cdot \phi_{j}dx\;\;{\rm strongly}\;{\rm in}\; L^{1}(0,T).
\end{equation}
In fact, we have
\begin{equation}\nonumber
\begin{split}
&\int_{\Omega_{f}}v_{\mu}\cdot\nabla d_{\mu}\cdot \phi_{j}dx-\int_{\Omega_{f}}v\cdot\nabla d\cdot \phi_{j}dx\\
=&\int_{\Omega_{f}}(v_{\mu}-v)\cdot\nabla d_{\mu}\phi_{j}dx+\int_{\Omega_{f}}v\cdot\nabla(d_{\mu}-d) \phi_{j}dx\\
\leq&\|v_{\mu}-v\|_{L^{3}(\Omega_{f})}\|\nabla d_{\mu}\|_{L^{2}(\Omega_{f})}\| \phi_{j}\|_{L^{6}(\Omega_{f})}+\int_{\Omega_{f}}v^{i}\phi^{k}_{j}\partial^{i}(d^{k}_{\mu}-d^{k}) dx\\
\leq&\|v_{\mu}-v\|_{H^{1-\epsilon}(\Omega_{f})}\|d_{\mu}\|_{H^{1}(\Omega_{f})}\| \phi_{j}\|_{H^{1}(\Omega_{f})}+\int_{\Omega_{f}}v^{i}\phi^{k}_{j}\partial^{i}(d^{k}_{\mu}-d^{k})dx\\
& \rightarrow 0, \quad \mu\rightarrow\infty,
\end{split}
\end{equation}
where we have used H\"{o}lder's inequality, $\eqref{1.12}_{1}$, $\eqref{2.16}_{2}$, $\eqref{2.16}_{3}$, $\eqref{2.17}_{1}$ and $\|v^{i}\phi^{k}_{j}\|_{L^{2}(\Omega_{f})}\leq C$. The $v^{i}$ denotes the $ith$ component of $v$ and $\phi^{k}_{j}$ is the $kth$ component of $\phi_{j}$. Similarly for \eqref{2.22}, we can easily get\\
\begin{equation}\label{2.23}
\int_{\Omega_{f}}v_{\mu}\cdot\nabla v_{\mu}\cdot \varphi_{j}dx\rightarrow\int_{\Omega_{f}}v\cdot\nabla v\cdot \varphi_{j}dx\;\;{\rm strongly }\;{\rm in}\;L^{1}(0,T).
 \end{equation}\\
Thirdly, as $\mu\rightarrow\infty$, we have
\begin{eqnarray}\label{2.24}
\begin{split}
&(\nabla v_{\mu},\nabla \varphi_{j})_{\Omega_{f}} \rightharpoonup (\nabla v,\nabla \varphi_{j})_{\Omega_{f}} \;{\rm weakly}  \; {\rm in}\;L^{2}(0, T),\\
&(\nabla d_{\mu},\nabla \phi_{j})_{\Omega} \rightharpoonup (\nabla d,\nabla \phi_{j})_{\Omega} \;{\rm weakly}  \; {\rm in}\;L^{2}(0, T),\\
&(\nabla \int_{0}^{t}v_{\mu}d\sigma,\nabla \varphi_{j})_{\Omega_{s}} \rightharpoonup (\nabla \int_{0}^{t}v d\sigma,\nabla \varphi_{j})_{\Omega_{s}} \;{\rm weakly*}\; {\rm in}\;L^{\infty}(0, T),\\
&(\nabla w_{0\mu},\nabla \varphi_{j})_{\Omega_{s}} \rightharpoonup (\nabla w_{0},\nabla \varphi_{j})_{\Omega_{s}} \;{\rm weakly*}\; {\rm in}\;L^{\infty}(0, T),\\
&(d_{\mu}e, \varphi_{j})_{\Omega_{f}}\rightharpoonup (d e, \varphi_{j})_{\Omega_{f}}\;{\rm weakly*}\;{\rm in
}\;L^{\infty}(0, T),\\
\end{split}
\end{eqnarray}
and
\begin{eqnarray}\label{2.25}
\begin{split}
&(v_{\mu,t}, \varphi_{j})_{\Omega}\rightarrow (v_{t}, \varphi_{j})_{\Omega}\;{\rm in
}\; \mathcal{D'}(0,T),\\
&(d_{\mu,t},  \phi_{j})_{\Omega}\rightarrow (d_{t}, \phi_{j})_{\Omega}\;{\rm in
}\; \mathcal{D'}(0,T).\\
\end{split}
\end{eqnarray}
Combining \eqref{2.19}, \eqref{2.21}, \eqref{2.22}, \eqref{2.23}, \eqref{2.24} and \eqref{2.25}, we get the following equations
\begin{eqnarray}\nonumber
\begin{split}
&(v_{t},\varphi_{j})_{\Omega}+\varepsilon a_{\Omega_{f}}(v,\varphi_{j})+ \mu a_{\Omega_{s}}(\int_{0}^{t}v d\sigma,\varphi_{j})+\mu a_{\Omega_{s}}(w_{0},\varphi_{j})\\
+&b_{\Omega_{f}}(v,v,\varphi_{j})-\frac{1}{2}\int_{\Gamma}(v\cdot n)v\cdot \varphi_{j}d\Gamma=(d e,\varphi_{j})_{\Omega_{f}}+(f,\varphi_{j})_{\Omega}\\
\end{split}
\end{eqnarray}
and
\begin{eqnarray}\nonumber
\begin{split}
&(d_{t},\phi_{j})_{\Omega}+k_{1} a_{\Omega_{f}}(d, \phi_{j})+ k_{2}a_{\Omega_{s}}(d, \phi_{j})+b_{\Omega_{f}}(v,d, \phi_{j})\\
&-\frac{1}{2}\int_{\Gamma}(v\cdot n)d\cdot \phi_{j}d\Gamma=(g, \phi_{j})_{\Omega}.\\
\end{split}
\end{eqnarray}
Taking into account the subsequence $\{\varphi_{j}\}_{j=1}^{\infty}$ and $\{\phi_{j}\}_{j=1}^{\infty}$ be an orthogonal basis of $V$ and $\Theta$ respectively, we deduce that
\begin{eqnarray}\nonumber
\begin{split}
(v_{t},\varphi)_{\Omega}&+\varepsilon a_{\Omega_{f}}(v,\varphi)+ \mu a_{\Omega_{s}}(\int_{0}^{t}v d\sigma,\varphi)+\mu a_{\Omega_{s}}(w_{0},\varphi)\\
&+b_{\Omega_{f}}(v,v,\varphi)-\frac{1}{2}\int_{\Gamma}(v\cdot n)v\cdot \varphi d\Gamma=(d e,\varphi)_{\Omega_{f}}+(f,\varphi)_{\Omega}\\
\end{split}
\end{eqnarray}
and
\begin{eqnarray}\nonumber
\begin{split}
(d_{t},\phi)_{\Omega}&+k_{1} a_{\Omega_{f}}(d,\phi)+ k_{2}a_{\Omega_{s}}(d,\phi)+b_{\Omega_{f}}(v,d,\phi)\\
&-\frac{1}{2}\int_{\Gamma} (v\cdot n)d\cdot \phi d\Gamma=(g,\phi)_{\Omega}\\
\end{split}
\end{eqnarray}
for all test functions $\varphi \in V$, $\phi \in \Theta$ satisfying $\|\varphi\|_{L^{2}(\Omega)}=\|\phi\|_{L^{2}(\Omega)}=1$.\\
3. Now, we need to check that
\begin{equation}\label{2.26}
v(0)=v_{0}(x)=\left\{\begin{aligned}
u_{0}(x),x\in \Omega_f\\
w_{1}(x), x\in \Omega_s
\end{aligned}\right.,\quad
d(0)=d_{0}(x)=\left\{\begin{aligned}
\rho_{0}(x), x\in \Omega_f\\
\theta_{0}(x), x\in \Omega_s
\end{aligned},\right.
\end{equation}
in the sense of $V$ and $\Theta$ respectively, and $w(0)=w_{0}(x)$ in the sense of $V_s$.

According to \eqref{2.15}, we know that
\begin{eqnarray}\nonumber
\begin{split}
&v_{\mu}\in C(0,T;V'),\,d_{\mu}\in C(0,T;\Theta'),\\
&\int_{0}^{t}v_{\mu} d\sigma\in C(0,T;V_{s}).\\
\end{split}
\end{eqnarray}
Then, we get that
\begin{eqnarray}\label{2.27}
\begin{split}
&v_{\mu}(0),\rightharpoonup v(0)\, weakly\, in\, V',\\
&d_{\mu}(0),\rightharpoonup v(0)\, weakly\, in\,\Theta',\\
&w_{\mu}(0),\rightharpoonup w(0)\, weakly\, in\, V_{s}.
\end{split}
\end{eqnarray}
On the other hand, \eqref{2.3} implies
\begin{eqnarray}\label{2.28}
\begin{split}
&v_{\mu}(0)\rightarrow v_{0},\, strongly\, in\, V,\\
&d_{\mu}(0)\rightarrow d_{0},\, strongly\, in\, \Theta,\\
&w_{\mu}(0)\rightarrow w_{0}, \, strongly\, in\, V_{s}.\\
\end{split}
\end{eqnarray}
We have proved \eqref{2.26} by \eqref{2.27} and \eqref{2.28}.\\
By combining the three steps mentioned earlier, passing to the limit as $m\rightarrow \infty$ in \eqref{2.2}, we have successfully completed the proof of the existence of global weak solutions. Meanwhile, as $m\rightarrow \infty$, we get the the energy inequality of weak solutions to the system $\eqref{1.3}$ that
\begin{equation}\nonumber
\begin{split}
&\|v\|_{L^{\infty}(0,T; {L^{2}(\Omega))}}^{2}+\|d\|_{L^{\infty}(0,T; {L^{2}(\Omega))}}^{2}+\|v\|_{L^{\infty}(0,T; {L^{2}(\Omega_{s}))}}^{2}+\mu\|\int_{0}^{t}v^{s} d\sigma+w_{0}\|_{L^{\infty}(0,T; V_{s})}^{2}\\
&+\varepsilon \|v\|_{L^{2}(0,T;V_{f})}^{2}+k_{1} \|d\|_{L^{2}(0,T;H^{1}(\Omega_{f}))}^{2}+k_{2}\|d\|_{L^{2}(0,T;H^{1}(\Omega_{s}))}^{2}\\
\leq &C(T)<\infty,\;t\in[0,T].\\
\end{split}
\end{equation}
Now, our focus shifts to establishing the uniqueness of the global weak solutions in the two-dimensional case.\\
\\
\textbf{Step 4: Uniqueness}

When d=2. Suppose that $(v, d)$, $(v_{*}, d_{*})$ be any two weak solutions to the Boussinesq-structure interaction system \eqref{1.3}. Let
$\chi=v-v_{*}$, $\psi=d-d_{*}$ and $\gamma(u,v,w)=\frac{1}{2}\int_{\Gamma}(u\cdot n)(vw)d\Gamma$.  We can derive the following equations:
\begin{equation}\label{2.29}
\begin{split}
(\chi_{t},\varphi)_{\Omega}&+\varepsilon a_{\Omega_{f}}(\chi,\varphi)+\mu a_{\Omega_{s}}(\int_{0}^{t}\chi d\sigma,\varphi)\\
&+b_{\Omega_{f}}(v,\chi,\varphi)+b_{\Omega_{f}}(\chi,v,\varphi)-b_{\Omega_{f}}(\chi,\chi,\varphi)\\
&-\gamma(v,\chi,\varphi)-\gamma(\chi,v,\varphi)+\gamma(\chi,\chi,\varphi)=(\psi e,\varphi)_{\Omega_{f}},
\end{split}
\end{equation}
and
\begin{equation}\label{2.30}
\begin{split}
(\psi_{t},\phi)_{\Omega}&+k_{1} a_{\Omega_{f}}(\psi,\phi)+k_{2} a_{\Omega_{s}}(\psi,\phi)\\
&+b_{\Omega_{f}}(v,\psi,\phi)+b_{\Omega_{f}}(\chi,d,\phi)-b_{\Omega_{f}}(\chi,\psi,\phi)\\
&-\gamma(v,\psi,\phi)-\gamma(\chi,d,\phi)+\gamma(\chi,\psi,\phi)=0,
\end{split}
\end{equation}
with the initial conditions
\begin{equation}\label{2.31}
(\chi, \psi)(0)=(0,0),
\end{equation}
where $t\in[0,T]$, $\varphi\in V$ and $ \phi\in \Theta$. Due to the presence of hyperbolic terms in the equations, we need to adopt a more difficult method by selecting a piecewise continuous function $\varsigma_{m}(t)$ on $[0,T]$ that satisfies the following conditions:
\begin{align}\nonumber
\begin{split}
\varsigma_{m}(t)= \left \{
\begin{array}{lr}
   1,      & t<s-\frac{2}{m},\\
   0,      & t>s-\frac{1}{m}.
\end{array}
\right.
\end{split}
\end{align}
Let $\varrho_{n}=\varrho_{n}(t)\in C_{0}^{\infty}$, such that $\varrho_{n}(t)=\varrho_{n}(-t)$, $\varrho_{n}(t)\geq0 $, and $\int_{-\infty}^{\infty}\varrho_{n}(t)dt=1$ when $\varrho(t)\in[-\frac{1}{n}, \frac{1}{n}]$.
When $n>2m$, we let

\begin{equation}\nonumber
\varphi(t)=((\varsigma_{m}\chi)*\varrho_{n}*\varrho_{n})\varsigma_{m},\,
\phi( t)=((\varsigma_{m}\psi)*\varrho_{n}*\varrho_{n})\varsigma_{m}, x\in \Omega
\end{equation}
where $*$ denotes the convolution to $t$. The $\chi$ and $\psi$ will extend to 0 outside of $[0,T]$.\\
Putting \eqref{2.29} and \eqref{2.30} together and integrating it from 0 to $T$, one has
\begin{equation}\label{2.32}
\sum_{i=1}^{13}I_{nm}^{i}=0,
\end{equation}
where
\begin{equation}\nonumber
\begin{split}
&I_{nm}^{1}=\int_{0}^{T}(\chi_{t},(\varsigma_{m}\chi)*\varrho_{n}*\varrho_{n})_{\Omega}\varsigma_{m}dt,\\
&I_{nm}^{2}=\int_{0}^{T}(\psi_{t},(\varsigma_{m}\psi)*\varrho_{n}*\varrho_{n})_{\Omega}\varsigma_{m}dt,\\
&I_{nm}^{3}=\varepsilon\int_{0}^{T}a_{\Omega_{f}}(\chi,(\varsigma_{m}\chi)*\varrho_{n}*\varrho_{n})\varsigma_{m}dt,\\
&I_{nm}^{4}=k_{1}\int_{0}^{T}a_{\Omega_{f}}(\psi,(\varsigma_{m}\psi)*\varrho_{n}*\varrho_{n})\varsigma_{m}dt,\\
&I_{nm}^{5}=\mu\int_{0}^{T}a_{\Omega_{s}}(\int_{0}^{t}\chi d\sigma,(\varsigma_{m}\chi)*\varrho_{n}*\varrho_{n})\varsigma_{m}dt,\\
&I_{nm}^{6}=k_{2}\int_{0}^{T}a_{\Omega_{s}}(\psi,(\varsigma_{m}\psi)*\varrho_{n}*\varrho_{n})\varsigma_{m}dt,\\
\end{split}
\end{equation}
and
\begin{equation}\nonumber
\begin{split}
&I_{nm}^{7}=\int_{0}^{T}b_{\Omega_{f}}(v,\chi,(\varsigma_{m}\chi)*\varrho_{n}*\varrho_{n})\varsigma_{m}dt-\int_{0}^{T}\gamma(v,\chi,(\varsigma_{m}\chi)*\varrho_{n}*\varrho_{n})\varsigma_{m}dt,\\
&I_{nm}^{8}=\int_{0}^{T}b_{\Omega_{f}}(v,\psi,(\varsigma_{m}\psi)*\varrho_{n}*\varrho_{n})\varsigma_{m}dt-\int_{0}^{T}\gamma(v,\psi,(\varsigma_{m}\psi)*\varrho_{n}*\varrho_{n})\varsigma_{m}dt,\\
&I_{nm}^{9}=\int_{0}^{T}b_{\Omega_{f}}(\chi,v,(\varsigma_{m}\chi)*\varrho_{n}*\varrho_{n})\varsigma_{m}dt-\int_{0}^{T}\gamma(\chi,v,(\varsigma_{m}\chi)*\varrho_{n}*\varrho_{n})\varsigma_{m}dt,\\
&I_{nm}^{10}=\int_{0}^{T}b_{\Omega_{f}}(\chi,d,(\varsigma_{m}\psi)*\varrho_{n}*\varrho_{n})\varsigma_{m}dt-\int_{0}^{T}\gamma(\chi,d,(\varsigma_{m}\psi)*\varrho_{n}*\varrho_{n})\varsigma_{m}dt,\\
&I_{nm}^{11}=\int_{0}^{T}b_{\Omega_{f}}(\chi,\chi,(\varsigma_{m}\chi)*\varrho_{n}*\varrho_{n})\varsigma_{m}dt-\int_{0}^{T}\gamma(\chi,\chi,(\varsigma_{m}\chi)*\varrho_{n}*\varrho_{n})\varsigma_{m}dt,\\
&I_{nm}^{12}=\int_{0}^{T}b_{\Omega_{f}}(\chi,\psi,(\varsigma_{m}\psi)*\varrho_{n}*\varrho_{n})\varsigma_{m}dt-\int_{0}^{T}\gamma(\chi,\psi,(\varsigma_{m}\psi)*\varrho_{n}*\varrho_{n})\varsigma_{m}dt,\\
&I_{nm}^{13}=-\int_{0}^{T}(\psi e,(\varsigma_{m}\chi)*\varrho_{n}*\varrho_{n})_{\Omega_{f}}\varsigma_{m}dt.\\
\end{split}
\end{equation}
Next, we will estimate each term on the left side of $\eqref{2.32}$ as follows.
By performing direct calculations, we can easily obtain that
\begin{equation}\label{2.33}
\begin{split}
I_{nm}^{1}=&\int_{0}^{T}((\varsigma_{m}\chi)_{t}-\varsigma_{m,t}\chi,(\varsigma_{m}\chi)*\varrho_{n}*\varrho_{n})_{\Omega}dt\\
=&\int_{0}^{T}((\varsigma_{m}\chi)_{t}*\varrho_{n},(\varsigma_{m}\chi)*\varrho_{n})_{\Omega}dt-\int_{0}^{T}\varsigma_{m,t}(\chi,(\varsigma_{m}\chi)*\varrho_{n}*\varrho_{n})_{\Omega}dt\\
=&-\int_{0}^{T}\varsigma_{m,t}(\chi,(\varsigma_{m}\chi)*\varrho_{n}*\varrho_{n})_{\Omega}dt.
\end{split}
\end{equation}
Using a same method as \eqref{2.33} for $I_{nm}^{2}$, $I_{nm}^{3}$ and $I_{nm}^{4}$. Letting $n\rightarrow \infty $, we have
\begin{equation}\nonumber
\lim_{n\rightarrow\infty}I_{nm}^{1}=-\int_{0}^{T}\varsigma_{m}\varsigma_{m,t}|\chi|_{\Omega}^{2}dt=I_{m}^{1},\,\lim_{n\rightarrow\infty}I_{nm}^{2}=-\int_{0}^{T}\varsigma_{m}\varsigma_{m,t}|\psi|_{\Omega}^{2}dt=I_{m}^{2}.
\end{equation}
Obviously, we can derive that
\begin{equation}\nonumber
\begin{split}
&\lim_{n\rightarrow\infty}I_{nm}^{3}=\varepsilon\int_{0}^{T}\varsigma_{m}^{2}a_{\Omega_{f}}(\chi,\chi)dt=I_{m}^{3},\\
&\lim_{n\rightarrow\infty}I_{nm}^{4}=k_{1}\int_{0}^{T}\varsigma_{m}^{2}a_{\Omega_{f}}(\psi,\psi)dt=I_{m}^{4},\\
&\lim_{n\rightarrow\infty}I_{nm}^{6}=k_{2}\int_{0}^{T}\varsigma_{m}^{2}a_{\Omega_{s}}(\psi,\psi)dt=I_{m}^{6},
\end{split}
\end{equation}
In order to estimate $I_{nm}^{5}$, letting
\begin{equation*}
\int_{0}^{t}\chi d\sigma=F(t),
\end{equation*}
then
\begin{equation*}
\begin{split}
I_{nm}^{5}=&\int_{0}^{T}\mu a_{\Omega_{s}}(\varsigma_{m}F, (\varsigma_{m}F_{t})*\varrho_{n}*\varrho_{n})dt\\
=&\int_{0}^{T}\mu a_{\Omega_{s}}((\varsigma_{m}F)*\varrho_{n},(\varsigma_{m}F)_{t}*\varrho_{n})dt-\int_{0}^{T}\mu a_{\Omega_{s}}(\varsigma_{m}F,(\varsigma_{m,t}F)*\varrho_{n}*\varrho_{n})dt\\
=&-\int_{0}^{T}\mu a_{\Omega_{s}}(\varsigma_{m}F,(\varsigma_{m,t}F)*\varrho_{n}*\varrho_{n})dt,
\end{split}
\end{equation*}
which can be deduced that
\begin{equation*}
\lim_{n\rightarrow\infty}I_{nm}^{5}=-\mu\int_{0}^{T}\varsigma_{m}\varsigma_{m,t}a_{\Omega_{s}}(F,F)dt=I_{m}^{5}.\\
\end{equation*}
By virtue of $b_{\Omega_{f}}(v,\chi,\chi)=\frac{1}{2}\int_{\Omega_{f}}v\cdot\nabla|\chi|^{2}dx=\frac{1}{2}\int_{\Gamma}|\chi|^{2}n\cdot vd\Gamma$, we have
\begin{equation}\nonumber
\begin{split}
&\lim_{n\rightarrow\infty}I_{nm}^{7}=\int_{0}^{T}[b_{\Omega_{f}}(v,\chi,\chi)-\gamma(v,\chi,\chi)]\varsigma_{m}^{2}dt=0,\\
&\lim_{n\rightarrow\infty}I_{nm}^{8}=\lim_{n\rightarrow\infty}I_{nm}^{11}=\lim_{n\rightarrow\infty}I_{nm}^{12}=0,\\
\end{split}
\end{equation}
and
\begin{equation}\nonumber
\begin{split}
&\lim_{n\rightarrow\infty}I_{nm}^{9}=\int_{0}^{T}[b_{\Omega_{f}}(\chi,v,\chi)-\gamma(\chi,v,\chi)]\varsigma_{m}^{2}dt=I_{m}^{9},\\
&\lim_{n\rightarrow\infty}I_{nm}^{10}=\int_{0}^{T}[b_{\Omega_{f}}(\chi,d,\psi)-\gamma(\chi,d,\psi)]\varsigma_{m}^{2}dt=I_{m}^{10},\\
&\lim_{n\rightarrow\infty}I_{nm}^{13}=-\int_{0}^{T}(\psi e,\chi)_{\Omega_{f}}\varsigma_{m}^{2}dt=I_{m}^{13}.
\end{split}
\end{equation}
Gathering those estimates, one checks
\begin{equation}\label{2.34}
I_{m}^{1}+I_{m}^{2}+I_{m}^{3}+I_{m}^{4}+I_{m}^{5}+I_{m}^{6}+I_{m}^{9}+I_{m}^{10}+I_{m}^{13}=0.
\end{equation}
Recalling Lebesgue Theorem, the following terms are held for $a.a$ $s\in [0,T]$
\begin{equation}\label{2.35}
\begin{split}
&-\int_{0}^{T}\varsigma_{m}\varsigma_{m,t}|\chi|_{\Omega}^{2}dt\rightarrow\frac{1}{2}\|\chi(s)\|_{L^{2}(\Omega)}^{2},\\
&-\int_{0}^{T}\varsigma_{m}\varsigma_{m,t}|\psi|_{\Omega}^{2}dt\rightarrow\frac{1}{2}\|\psi(s)\|_{L^{2}(\Omega)}^{2},\\
&-\mu\int_{0}^{T}\varsigma_{m}\varsigma_{m,t}a_{\Omega_{s}}(F ,F)dt\rightarrow\frac{\mu}{2}\|\nabla F(s)\|_{L^{2}(\Omega_{s})}^{2}.\\
\end{split}
\end{equation}
By the fact
\begin{equation}\label{2.36}
\begin{split}
&a_{\Omega_{f}}(\chi,\chi)=\|\chi\|_{H^{1}(\Omega_{f})}^{2}-\|\chi\|_{L^{2}(\Omega_{f})}^{2},\\
&a_{\Omega}(\psi,\psi)=\|\psi\|_{H^{1}(\Omega)}^{2}-\|\psi\|_{L^{2}(\Omega)}^{2}.
\end{split}
\end{equation}
Combining \eqref{2.34}, \eqref{2.35} and \eqref{2.36}, then we have
\begin{equation}\label{2.37}
\begin{split}
\frac{1}{2}&(\|\chi(s)\|_{L^{2}(\Omega)}^{2}+\|\psi(s)\|_{L^{2}(\Omega)}^{2}+\mu\|\nabla F(s)\|_{L^{2}(\Omega_{s})}^{2})\\
&+\int_{0}^{s}[\varepsilon(\|\chi(t)\|_{H^{1}(\Omega_{f})}^{2}-\|\chi(t)\|_{L^{2}(\Omega_{f})}^{2})+k_{1}(\|\psi(t)\|_{H^{1}(\Omega_{f})}^{2}-\|\psi(t)\|_{L^{2}(\Omega_{f})}^{2})\\
&+k_{2}(\|\psi(t)\|_{H^{1}(\Omega_{s})}^{2}-\|\psi(t)\|_{L^{2}(\Omega_{s})}^{2})]dt-\int_{0}^{s}(\psi e, \chi)_{\Omega_{f}}dt\\
&+\int_{0}^{s}[b_{\Omega_{f}}(\chi, v, \chi)-\gamma(\chi, v, \chi)]+[b_{\Omega_{s}}(\chi, d, \psi)-\gamma(\chi, d, \psi)]dt=0.
\end{split}
\end{equation}
Furthermore, we can write that
\begin{equation}\label{2.38}
\begin{split}
&\frac{1}{2}(\|\chi(s)\|_{L^{2}(\Omega)}^{2}+\|\psi(s)\|_{L^{2}(\Omega)}^{2}+\mu\|\nabla F(s)\|_{L^{2}(\Omega_{s})}^{2})\\
&+\int_{0}^{s}\varepsilon\|\chi(t)\|_{H^{1}(\Omega_{f})}^{2}+k_{1}\|\psi(t)\|_{H^{1}(\Omega_{f})}^{2}+k_{2}\|\psi(t)\|_{H^{1}(\Omega_{s})}^{2}dt\\
\leq&\int_{0}^{s}(\varepsilon\|\chi(t)\|_{L^{2}(\Omega_{f})}^{2}+k_{1}\|\psi(t)\|_{L^{2}(\Omega_{f})}^{2}+k_{2}\|\psi(t)\|_{L^{2}(\Omega_{s})}^{2})dt+|\int_{0}^{s}(\psi e, \chi)_{\Omega_{f}} dt|\\
&+|\int_{0}^{s}[b_{\Omega_{f}}(\chi, v, \chi)-\gamma(\chi, v, \chi)]dt|+|\int_{0}^{s}[b_{\Omega_{s}}(\chi, d, \psi)-\gamma(\chi,d, \psi)]dt|\\
\leq&\int_{0}^{s}(\varepsilon\|\chi(t)\|_{L^{2}(\Omega_{f})}^{2}+k_{1}\|\psi(t)\|_{L^{2}(\Omega_{f})}^{2}+k_{2}\|\psi(t)\|_{L^{2}(\Omega_{s})}^{2})dt+|\int_{0}^{s}(\psi e, \chi)_{\Omega_{f}} dt|\\
&+\int_{0}^{s}|b_{\Omega_{f}}(\chi, v, \chi)|+|b_{\Omega_{s}}(\chi, d, \psi)|dt+\int_{0}^{s}|\gamma(\chi, v, \chi)|+|\gamma(\chi, d, \psi)|dt\\
=&\int_{0}^{s}\sum_{i=1}^{6} J_{i}dt.
\end{split}
\end{equation}
Next, we will start estimating from the second term in \eqref{2.38}.
\begin{equation}\label{2.39}
J_{2}=\int_{\Omega_{f}}\psi e\cdot \chi dt\leq \frac{1}{2}\|\psi\|_{L^{2}(\Omega_{f})}^{2}+ \frac{1}{2}\|\chi\|_{L^{2}(\Omega_{f})}^{2},
\end{equation}
and
\begin{equation}\label{2.40}
\begin{split}
J_{3}&=|b_{\Omega_{f}}(\chi, v, \chi)|\\
&\leq C\|\chi\|^{2}_{L^{4}(\Omega_{f})}\|\nabla v\|_{L^{2}(\Omega_{f})}\\
&\leq C\|\chi\|_{L^{2}(\Omega_{f})}\|\nabla \chi\|_{L^{2}(\Omega_{f})}\|v\|_{H^{1}(\Omega_{f})}\\
&\leq C\|\chi\|_{L^{2}(\Omega_{f})}\|\chi\|_{H^{1}(\Omega_{f})}\|v\|_{H^{1}(\Omega_{f})}\\
&\leq \frac{\varepsilon}{8}\|\chi\|_{H^{1}(\Omega_{f})}^{2}+C\|\chi\|_{L^{2}(\Omega_{f})}^{2}\|v\|_{H^{1}(\Omega_{f})}^{2},\\
\end{split}
\end{equation}
and
\begin{equation}\label{2.41}
\begin{split}
J_{4}&=|b_{\Omega_{f}}(\chi, d, \psi)|\\
&\leq C\|\chi\|_{L^{4}(\Omega_{f})}\|\nabla d\|_{L^{2}(\Omega_{f})}\|\psi\|_{L^{4}(\Omega_{f})}\\
&\leq C\|\chi\|_{L^{2}(\Omega_{f})}^{1/2}\|\chi\|_{H^{1}(\Omega_{f})}^{1/2}\|\psi\|_{L^{2}(\Omega_{f})}^{1/2}\|\psi\|_{H^{1}(\Omega_{f})}^{1/2} \|d\|_{H^{1}(\Omega_{f})}\\
&\leq \frac{k_{1}}{4}\| \psi\|_{H^{1}(\Omega_{f})}^{2}+C\|\chi\|_{L^{2}(\Omega_{f})}^{2/3}\|\chi\|_{H^{1}(\Omega_{f})}^{2/3}\|\psi\|_{L^{2}(\Omega_{f})}^{2/3} \|d\|^{4/3}_{H^{1}(\Omega_{f})}\\
&\leq \frac{k_{1}}{4}\| \psi\|_{H^{1}(\Omega_{f})}^{2}+\frac{\varepsilon}{8}\|\chi\|_{H^{1}(\Omega_{f})}^{2}+C(\|\chi\|_{L^{2}(\Omega_{f})}^{2}+\|\psi\|_{L^{2}(\Omega_{f})}^{2})\|d\|_{H^{1}(\Omega_{f})}^{2},\\
\end{split}
\end{equation}
where we have used H\"{o}lder's inequality, Young's inequality and $\eqref{1.12}$ with $\mathrm{d}=2$.
\begin{equation}\label{2.42}
\begin{split}
J_{5}&=|\gamma(\chi, v, \chi)|\leq C\|\chi|_{\Gamma}\|_{L^{3}(\Gamma)}^{2}\|v|_{\Gamma}\|_{L^{3}(\Gamma)}\\
&\leq C\|\chi\|_{H^{1}(\Omega_{f})}^{4/3}\|\chi\|_{L^{2}(\Omega_{f})}^{2/3}\|v|_{\Gamma}\|_{L^{3}(\Gamma)}\\
&\leq\frac{\varepsilon}{8}\|\chi\|_{H^{1}(\Omega_{f})}^{2}+C\|\chi\|_{L^{2}(\Omega_{f})}^{2}\|v|_{\Gamma}\|_{L^{3}(\Gamma)}^{3},
\end{split}
\end{equation}
and
\begin{equation}\label{2.43}
\begin{split}
J_{6}&=|\gamma(\chi, d, \psi)|\leq C\|\chi|_{\Gamma}\|_{L^{3}(\Gamma)} \|\psi|_{\Gamma}\|_{L^{3}(\Gamma)}\|d|_{\Gamma}\|_{L^{3}(\Gamma)}\\
&\leq C\|\chi\|_{H^{1}(\Omega_{f})}^{2/3}\|\chi\|_{L^{2}(\Omega_{f})}^{1/3}\|\psi\|_{H^{1}(\Omega_{f})}^{2/3}\|\psi\|_{L^{2}(\Omega_{f})}^{1/3}\|d|_{\Gamma}\|_{L^{3}(\Gamma)}\\
&\leq\frac{\varepsilon}{8}\|\chi\|_{H^{1}(\Omega_{f})}^{2}+C\|\chi\|_{L^{2}(\Omega_{f})}^{1/2}\|\psi\|_{H^{1}(\Omega_{f})}\|\psi\|_{L^{2}(\Omega_{f})}^{1/2}\|d|_{\Gamma}\|_{L^{3}(\Gamma)}^{3/2}\\
&\leq\frac{\varepsilon}{8}\|\chi\|_{H^{1}(\Omega_{f})}^{2}+\frac{k_{1}}{4}\|\psi\|_{H^{1}(\Omega_{f})}^{2}+C(\|\chi\|_{L^{2}(\Omega_{f})}^{2}+\|\psi\|_{L^{2}(\Omega_{f})}^{2})\|d|_{\Gamma}\|_{L^{3}(\Gamma)}^{3},
\end{split}
\end{equation}
where we have used H\"{o}lder's inequality, Young's inequality and $\eqref{1.11}_{2}$ with $\mathrm{d}=2$.
Putting \eqref{2.39}-\eqref{2.43} into \eqref{2.38}, we obtain that
\begin{equation}\label{2.44}
\begin{split}
&(\|\chi\|_{L^{2}(\Omega)}^{2}+\|\psi\|_{L^{2}(\Omega)}^{2}+\mu\|\nabla F(s)\|_{L^{2}(\Omega_{s})}^{2})+\int_{0}^{s}(\varepsilon\|\chi\|_{H^{1}(\Omega_{f})}^{2}+k_{1}\|\psi\|_{H^{1}(\Omega_{f})}^{2}+k_{2}\|\psi\|_{H^{1}(\Omega_{s})}^{2})dt\\
\leq&2\int_{0}^{s}(\varepsilon\|\chi\|_{L^{2}(\Omega_{f})}^{2}+k_{1}\|\psi\|_{L^{2}(\Omega_{f})}^{2}+k_{2}\|\psi\|_{L^{2}(\Omega_{s})}^{2})dt+\int_{0}^{s}\|\psi\|_{L^{2}(\Omega_{f})}^{2}+ \|\chi\|_{L^{2}(\Omega_{f})}^{2}dt\\
&+\int_{0}^{s}C(\|d\|_{H^{1}(\Omega_{f})}^{2}+\|d|_{\Gamma}\|_{L^{3}(\Gamma)}^{3}+\|v\|_{H^{1}(\Omega_{f})}^{2}+\|v|_{\Gamma}\|_{L^{3}(\Gamma)}^{3})(\|\chi\|_{L^{2}(\Omega_{f})}^{2}+\|\psi\|_{L^{2}(\Omega_{f})}^{2})dt\\
\leq&\int_{0}^{s}CM(t)(\|\chi\|_{L^{2}(\Omega)}^{2}+\|\psi\|_{L^{2}(\Omega)}^{2})dt,\\
\end{split}
\end{equation}
where
\begin{equation}\nonumber
\begin{split}
&M(t)=C'+\|d\|_{H^{1}(\Omega_{f})}^{2}+\|d|_{\Gamma}\|_{L^{3}(\Gamma)}^{3}+\|v\|_{H^{1}(\Omega_{f})}^{2}+\|v|_{\Gamma}\|_{L^{3}(\Gamma)}^{3},\\
&C'=max\{2\varepsilon, 2k_{1},2k_{2},1\}.
\end{split}
\end{equation}
However, one obtains that $M(t)\in L^{1}(0,T)$ and
\begin{equation}\label{2.45}
\|\chi\|_{L^{2}(\Omega)}^{2}+\|\psi\|_{L^{2}(\Omega)}^{2}\leq \int_{0}^{s}CM(t)(\|\chi\|_{L^{2}(\Omega)}^{2}+\|\psi\|_{L^{2}(\Omega)}^{2})dt.
\end{equation}
Applying Gronwall's inequality to \eqref{2.45} and employing the initial conditions \eqref{2.31}, we get
\begin{equation}\label{2.46}
\chi(t)=0,\; \psi(t)=0.
\end{equation}
By combining the four steps outlined above, we have successfully proven the existence of global weak solutions in two or three-dimensional space. Additionally, we have established the uniqueness of these solutions specifically in the two-dimensional case.

\section{\bf Proof of Theorem 1.2}

We have proved existence of weak solutions to the system \eqref{1.1} in Theorem \ref{thm 1.1.}, the rest is to obtain the more regularity estimates in two-dimensional space for the approximation solution $\{v_{m},d_{m}\}$ given by \eqref{2.1}. Hence, we continue take $\{\varphi_{j},\phi_{j}\}$ as the orthogonal basis sequence in space $(H^{2}(V)\times H^{2}(\Theta))$. The initial conditions \eqref{2.3} satisfies
\begin{equation*}
\begin{cases}
v_{m}(0)\rightarrow u_{0} \; strongly \;in \;V_{f}\cap H^{2}(\Omega_{f}),
d_{m}(0)\rightarrow \rho_{0}\; strongly \;in \;H^{2}(\Omega_{f}),\\
w_{m}(0)\rightarrow w_{0}\; strongly \;in \; V_{s}\cap H^{2}(\Omega_{s}),\\
v_{m}(0)\rightarrow w_{1}\; strongly \;in\;V_{s}, d_{m}(0)\rightarrow \theta_{0}\; strongly \;in \; H^{2}(\Omega_{s})\\
\end{cases}
\end{equation*}
and
\begin{equation}\label{3.1}
\begin{cases}
 \mu\frac{\partial w_{0}}{\partial n}-\varepsilon\frac{\partial u_{0}}{\partial n}+\frac{1}{2}(u_{0}\cdot n)u_{0}=(\mu\frac{\partial w_{0}}{\partial n}-\varepsilon\frac{\partial u_{0}}{\partial n}+\frac{1}{2}(u_{0}\cdot n)\cdot u_{0})\cdot n, x\in\Gamma,\\
 k_{2}\frac{\partial \theta_{0}}{\partial n}-k_{1}\frac{\partial \rho_{0}}{\partial n}+\frac{1}{2}(u_{0}\cdot n)\rho_{0}=0,x\in\Gamma.
\end{cases}
\end{equation}
\textbf{Step 1: A prior estimates of the initial data $(v_{m,t}(0),d_{m,t}(0))$.}\\

The proposal of compatibility conditions helps us solve nonlinear boundary integral terms. Integrating by parts, using H\"{o}lder's inequality and the trace property to $\eqref{2.2}_{1}$, and together with $\eqref{3.1}_{1}$, we can get the uniform estimates for $v_{m,t}(0)$ on $m$.
\begin{equation}\label{3.2}
\begin{split}
&\|v_{m,t}(0)\|^{2}_{L^{2}(\Omega)}\\
=&\varepsilon \int_{\Omega_{f}}\Delta v_{m}(0)\cdot v_{m,t}(0)dx+\mu \int_{\Omega_{s}}\Delta w_{0}\cdot v_{m,t}(0)dx+\int_{\Omega_{f}}d_{m}(0)e\cdot v_{m,t}(0)dx\\
&-\int_{\Gamma}(\varepsilon\frac{\partial v_{m}(0)}{\partial n}-\mu\frac{\partial w_{0}}{\partial n})\cdot v_{m,t}(0)d\Gamma -b_{\Omega_{f}}(v_{m}(0),v_{m}(0),v_{m,t}(0))\\
&+\frac{1}{2}\int_{\Gamma}(v_{m}(0)\cdot n)v_{m}(0)\cdot v_{m,t} (0) d\Gamma+\int _{\Omega}f(0)\cdot v_{m,t}(0)dx\\
\leq&\varepsilon\|\Delta v_{m}(0)\|_{L^{2}(\Omega_{f})}\|v_{m,t}(0)\|_{L^{2}(\Omega_{f})}+\mu\|\Delta w_{0}\|_{L^{2}(\Omega_{s})}\|v_{m,t}(0)\|_{L^{2}(\Omega_{s})}\\
&+\|d_{m}(0)\|_{L^{2}(\Omega_{f})}\|v_{m,t}(0)\|_{L^{2}(\Omega_{f})}+\|((\varepsilon\frac{\partial v_{m}(0)}{\partial n}-\mu\frac{\partial w_{0}}{\partial n})\cdot n)|_{\Gamma}\|_{H^{\tfrac{1}{2}}(\Gamma)}\|(n\cdot v_{m,t}(0))|_{\Gamma}\|_{H^{-\tfrac{1}{2}}(\Gamma)}\\
&+\|v_{m}(0)\cdot\nabla v_{m}(0)\|_{L^{2}(\Omega_{f})}\|v_{m,t}(0)\|_{L^{2}(\Omega_{f})}+\|v_{m}(0)|_{\Gamma}|^{2}\|_{H^{\tfrac{1}{2}}(\Gamma)}\|n\cdot v_{m,t}(0)|_{\Gamma}\|_{H^{-\tfrac{1}{2}}(\Gamma)}\\
&+\|f(0)\|_{L^{2}(\Omega)}\|v_{m,t}(0)\|_{L^{2}(\Omega)}\\
\leq&C(\|(\Delta v_{m}(0),d_{m}(0))\|_{L^{2}(\Omega_{f})}+\|\Delta w_{0} \|_{L^{2}(\Omega_{s})}+\|f(0)\|_{L^{2}(\Omega)}    )\times\|v_{m,t}(0)\|_{L^{2}(\Omega)}\\
&+(\|\frac{\partial v_{m}(0)}{\partial n}\|_{H^{1}(\Omega_{f})}+\|\frac{\partial w_{0}}{\partial n}\|_{H^{1}(\Omega_{s})})\|v_{m,t}(0)\|_{L^{2}(\Omega_{f})}+\|v_{m}(0)\|^{2}_{H^{2}(\Omega_{f})}\|v_{m,t}(0)\|_{L^{2}(\Omega_{f})}\\
\leq&C(\|u_{0}\|_{H^{2}(\Omega_{f})}+\|u_{0}\|^{2}_{H^{2}(\Omega_{f})}+\| w_{0}\|_{H^{2}(\Omega_{s})}+\|\rho_{0}\|_{L^{2}(\Omega_{f})}+\|f(0)\|_{L^{2}(\Omega)})\times\|v_{m,t}(0)\|_{L^{2}(\Omega)},
\end{split}
\end{equation}
which implies that
\begin{equation}\label{3.3}
\begin{split}
&\|v_{m,t}|_{t=0}\|_{L^{2}(\Omega)}\\
\leq&C(\|u_{0}\|_{H^{2}(\Omega_{f})}+\|u_{0}\|^{2}_{H^{2}(\Omega_{f})}+\| w_{0}\|_{H^{2}(\Omega_{s})}+\|\rho_{0}\|_{L^{2}(\Omega_{f})}+\|f(0)\|_{L^{2}(\Omega)})\\
\leq &C<\infty.\\
\end{split}
\end{equation}
This means that, by $\eqref{1.3}_{1}$ and $\eqref{1.3}_{4}$ at $t=0$,
\begin{equation}\label{3.4}
\begin{split}
&\|v_{m,t}|_{t=0}\|_{L^{2}(\Omega)}\\
\leq& C(\|\varepsilon\Delta u_{0}+\rho_{0} e-u_{0}\cdot \nabla u_{0}-\nabla p(t=0)\|_{L^{2}(\Omega_{f})}+\|\mu\Delta w_{0}\|_{L^{2}(\Omega_{s})}+\|f(0)\|_{L^{2}(\Omega)})\\
\leq&C<\infty,
\end{split}
\end{equation}
where the pressure $p_{0}\in H^{1}(\Omega_{f})$ is given that
\begin{equation}\label{3.5}
\begin{cases}
 -\Delta p_{0}={\rm div}(u_{0}\cdot\nabla u_{0}-\rho_{0} e-f_{1}(0)), x\in \Omega_{f},\\
 p_{0}=(\varepsilon\frac{\partial u_{0}}{\partial n}-\mu\frac{\partial w_{0}}{\partial n}-\frac{1}{2}(u_{0}\cdot n)u_{0}) n, x\in\Gamma.
\end{cases}
\end{equation}
By combining the boundary condition $\eqref{1.3}_{7}$ and $\eqref{3.5}$, we know that the compatibility condition $\eqref{1.14}_{1}$ is necessary.
Using a similar method to $\eqref{2.1}_{2}$ and $\eqref{3.1}_{2}$, and combining $\eqref{1.3}_{3}$ and $\eqref{1.3}_{5}$, we can obtain the uniform estimates for $d_{m,t}(0)$ on $m$.
\begin{equation}\label{3.6}
\begin{split}
&\|d_{m,t}|_{t=0}\|^{2}_{L^{2}(\Omega)}\\
\leq& C(\|k_{1}\Delta\rho_{0}-u_{0}\cdot \nabla \rho_{0}\|_{L^{2}(\Omega_{f})}+\|k_{2}\Delta \theta_{0}\|_{L^{2}(\Omega_{s})}+\|g(0)\|_{L^{2}(\Omega)})\leq \infty.
\end{split}
\end{equation}\\
Furthermore, we can prove that the compatibility condition $\eqref{1.14}_{2}$ is required.\\
\\
\textbf{Step 2: A prior estimates of  $(v_{m,t}(t),d_{m,t}(t))$.}\\

Differentiating the equality $\eqref{2.2}_{1}$ with respect to $t$ and multiplying $\partial_{t}a^{m}_{j}(t)$, then summing from $j=1$ to $m$, we have
\begin{equation}\label{3.7}
\begin{split}
\frac{1}{2}\frac{d}{dt}&(\|v_{m,t}\|_{L^{2}(\Omega)}^{2}+\mu\|\nabla v_{m}\|_{L^{2}(\Omega_{s})}^{2})+\varepsilon \|\nabla v_{m,t}\|_{L^{2}(\Omega_{f})}^{2}\\
=&-\int_{\Omega_{f}}(v_{m}^{j}\partial^{j}v_{m}^{i})_{t}v^{i}_{m,t}dx+\frac{1}{2}\int_{\Gamma}(v^{j}_{m}v^{i}_{m})_{t}n^{j} v^{i}_{m,t}d\Gamma\\
&+\int_{\Omega_{f}}d_{m,t}e\cdot v_{m,t}dx+\int_{\Omega}f_{t}\cdot v_{m,t}dx.\\
\end{split}
\end{equation}
Recalling that $\int_{\Omega_{f}}v_{m}^{j}\partial^{j}v_{m,t}^{i}v^{i}_{m,t}dx=\frac{1}{2}\int_{\Gamma}v^{j}_{m}v^{i}_{m,t}n^{j} v^{i}_{m,t}d\Gamma$, from \eqref{3.7}, we obtain
\begin{equation}\label{3.8}
\begin{split}
\frac{1}{2}\frac{d}{dt}&(\|v_{m,t}\|_{L^{2}(\Omega)}^{2}+\mu\|\nabla v_{m}\|_{L^{2}(\Omega_{s})}^{2})+\varepsilon \|\nabla v_{m,t}\|_{L^{2}(\Omega_{f})}^{2}\\
=&-\int_{\Omega_{f}}v_{m,t}^{j}\partial^{j}v_{m}^{i}v^{i}_{m,t}dx+\frac{1}{2}\int_{\Gamma}v^{j}_{m,t}v^{i}_{m}n^{j} v^{i}_{m,t}d\Gamma\\
&+\int_{\Omega_{f}}d_{m,t}e\cdot v_{m,t}dx+\int_{\Omega}f_{t}\cdot v_{m,t}dx\\
\triangleq& \sum_{i=1}^{4}I_{i}.\\
\end{split}
\end{equation}
Next, we will estimate each term on the right-hand side of equation \eqref{3.8}.
\begin{equation}\label{3.9}
\begin{split}
I_{1}&\leq \|v_{m,t}\|^{2}_{L^{4}(\Omega_{f})}\|\nabla v_{m}\|_{L^{2}(\Omega_{f})}\\
&\leq C\|v_{m,t}\|_{L^{2}(\Omega_{f})}\|v_{m,t}\|_{H^{1}(\Omega_{f})}\|\nabla v_{m}\|_{L^{2}(\Omega_{f})}\\
&\leq\frac{\varepsilon}{8}\| v_{m,t}\|_{H^{1}(\Omega_{f})}^{2}+C\|v_{m,t}\|_{L^{2}(\Omega_{f})}^{2}\|\nabla v_{m}\|_{L^{2}(\Omega_{f})}^{2},\\
\end{split}
\end{equation}
and
\begin{equation}\label{3.10}
\begin{split}
I_{2}&\leq \|v_{m,t}|_{\Gamma}\|_{L^{3}(\Gamma)}^{2}\|v_{m}|_{\Gamma}\|_{L^{3}(\Gamma)}\\
&\leq C \|v_{m,t}|_{\Gamma}\|_{H^{1/6}(\Gamma)}^{2}\|v_{m}|_{\Gamma}\|_{H^{1/6}(\Gamma)}\\
&\leq C\|v_{m,t}\|_{H^{2/3}(\Omega_{f})}^{2}\|v\|_{H^{2/3}(\Omega_{f})}\\
&\leq C \|v_{m,t}\|_{L^{2}(\Omega_{f})}^{2/3}\|v_{m,t}\|_{H^{1}(\Omega_{f})}^{4/3} \|v_{m}\|_{L^{2}(\Omega_{f})}^{1/3}\|v_{m}\|_{H^{1}(\Omega_{f})}^{2/3}\\
&\leq\frac{\varepsilon}{8}\|v_{m,t}\|_{H^{1}(\Omega_{f})}^{2}+C\|v_{m,t}\|_{L^{2}(\Omega_{f})}^{2} \|v_{m}\|_{L^{2}(\Omega_{f})}\|v_{m}\|_{H^{1}(\Omega_{f})}^{2}\\
&\leq\frac{\varepsilon}{8}\|v_{m,t}\|_{H^{1}(\Omega_{f})}^{2}+C(T)\|v_{m,t}\|_{L^{2}(\Omega_{f})}^{2} \|\nabla v_{m}\|_{L^{2}(\Omega_{f})}^{2},\\
\end{split}
\end{equation}
where we have used H\"{o}lder's inequality, interpolation inequality, trace theorem, Young's inequality and \eqref{2.4}. A similar estimate can be applied as follows:
\begin{equation}\label{3.11}
I_{3}+I_{4}\leq \frac{1}{2}\|d_{m,t}\|_{L^{2}(\Omega_{f})}^{2}+\|v_{m,t}\|_{L^{2}(\Omega)}^{2}+\frac{1}{2}\|f_{t}\|_{L^{2}(\Omega)}^{2}.
\end{equation}
Differentiating the equality $\eqref{2.2}_{2}$ with respect to $t$ and multiplying the results with $\partial_{t}b^{m}_{j}(t)$, then summing from $j=1$ to $m$, we deduce that
\begin{equation}\label{3.12}
\begin{split}
&\frac{1}{2}\frac{d}{dt}\|d_{m,t}\|_{L^{2}(\Omega)}^{2}+k_{1}\|\nabla d_{m,t}\|_{L^{2}(\Omega_{f})}^{2}+k_{2} \|\nabla d_{m,t}\|_{L^{2}(\Omega_{s})}^{2}\\
=&-\int_{\Omega_{f}}(v_{m}^{j}\partial^{j}d_{m}^{i})_{t}\cdot d^{i}_{m,t}dx+\frac{1}{2}\int_{\Gamma}(v_{m}^{j}d_{m}^{i})_{t} n^{j}d^{i}_{m,t}d\Gamma+\int_{\Omega}g_{t}\cdot d_{m,t}dx.\\
\end{split}
\end{equation}
Recalling that $\int_{\Omega_{f}}v_{m}^{j}\partial^{j}d_{m,t}^{i}d^{i}_{m,t}dx=\frac{1}{2}\int_{\Gamma}v^{j}_{m}d^{i}_{m,t}n^{j} d^{i}_{m,t}d\Gamma$, from \eqref{3.12}, we obtain

\begin{equation}\label{3.13}
\begin{split}
&\frac{1}{2}\frac{d}{dt}\|d_{m,t}\|_{L^{2}(\Omega)}^{2}+k_{1}\|\nabla d_{m,t}\|_{L^{2}(\Omega_{f})}^{2}+k_{2} \|\nabla d_{m,t}\|_{L^{2}(\Omega_{s})}^{2}\\
=&-\int_{\Omega_{f}}v_{m,t}^{j}\partial^{j}d_{m}^{i} d^{i}_{m,t}dx+\frac{1}{2}\int_{\Gamma}v_{m,t}^{j}d_{m}^{i} n^{j}d^{i}_{m,t}d\Gamma+\int_{\Omega}g_{t}\cdot d_{m,t}dx\\
\triangleq& \sum_{i=1}^{3}J_{i}.\\
\end{split}
\end{equation}
We need to estimate every term in \eqref{3.13}.
\begin{equation}\label{3.14}
\begin{split}
J_{1}\leq &\|v_{m,t}\|_{L^{4}(\Omega_{f})}\|d_{m,t}\|_{L^{4}(\Omega_{f})}\|\nabla d_{m}\|_{L^{2}(\Omega_{f})}\\
\leq & C\|v_{m,t}\|_{L^{2}(\Omega_{f})}^{1/2}\|v_{m,t}\|_{H^{1}(\Omega_{f})}^{1/2}\|d_{m,t}\|_{L^{2}(\Omega_{f})}^{1/2}\|d_{m,t}\|_{H^{1}(\Omega_{f})}^{1/2}\| \nabla d_{m}\|_{L^{2}(\Omega_{f})}\\
\leq & \frac{\varepsilon}{8}\|v_{m,t}\|_{H^{1}(\Omega_{f})}^{2}+\frac{k_{1}}{4}\|d_{m,t}\|_{H^{1}(\Omega_{f})}^{2}+C\|v_{m,t}\|_{L^{2}(\Omega_{f})}\|d_{m,t}\|_{L^{2}(\Omega_{f})}\|\nabla d_{m}\|^{2}_{L^{2}(\Omega_{f})},\\
\leq &\frac{\varepsilon}{8}\| v_{m,t}\|_{H^{1}(\Omega_{f})}^{2}+\frac{k_{1}}{4}\|d_{m,t}\|_{H^{1}(\Omega_{f})}^{2}+
C(\|v_{m,t}\|_{L^{2}(\Omega_{f})}^{2}
+\|d_{m,t}\|_{L^{2}(\Omega_{f})}^{2})\|\nabla d_{m}\|^{2}_{L^{2}(\Omega_{f})},\\
\end{split}
\end{equation}
and
\begin{equation}\label{3.15}
\begin{split}
J_{2}\leq &\|d_{m,t}|_{\Gamma}\|_{L^{3}(\Gamma)}\|v_{m,t}|_{\Gamma}\|_{L^{3}(\Gamma)} \|d_{m}|_{\Gamma}\|_{L^{3}(\Gamma)}\\
\leq&C\|d_{m,t}\|_{H^{2/3}(\Omega_{f})}\|v_{m,t}\|_{H^{2/3}(\Omega_{f})}\|v_{m}\|_{H^{2/3}(\Omega_{f})}\\
\leq&C \|d_{m,t}\|_{L^{2}(\Omega_{f})}^{1/3}\|d_{m,t}\|_{H^{1}(\Omega_{f})}^{2/3}\|v_{m,t}\|_{L^{2}(\Omega_{f})}^{1/3}\|v_{m,t}\|_{H^{1}(\Omega_{f})}^{2/3} \|d_{m}\|_{L^{2}(\Omega_{f})}^{1/3}\|d_{m}\|_{H^{1}(\Omega_{f})}^{2/3}\\
\leq &\frac{k_{1}}{4}\|d_{m,t}\|_{H^{1}(\Omega_{f})}^{2}+\frac{\varepsilon}{8}\|v_{m,t}\|_{H^{1}(\Omega_{f})}^{2}+C\|d_{m,t}\|_{L^{2}(\Omega_{f})}\|v_{m,t}\|_{L^{2}(\Omega_{f})}\|d_{m}\|_{L^{2}(\Omega_{f})}\|d_{m}\|^{2}_{H^{1}(\Omega_{f})},\\
\leq &\frac{k_{1}}{4}\|d_{m,t}\|_{H^{1}(\Omega_{f})}^{2}+\frac{\varepsilon}{8}\|v_{m,t}\|_{H^{1}(\Omega_{f})}^{2}+C(T)(\|d_{m,t}\|^{2}_{L^{2}(\Omega_{f})}+\|v_{m,t}\|^{2}_{L^{2}(\Omega_{f})})\|\nabla d_{m}\|^{2}_{L^{2}(\Omega_{f})},\\
\end{split}
\end{equation}
where we have utilized the same methods as \eqref{3.9}, \eqref{3.10}. And it is easy to obtain the following results
that
\begin{equation}\label{3.16}
 J_{3}\leq \frac{1}{2}\|d_{m,t}\|_{L^{2}(\Omega)}^{2}+\frac{1}{2}\|g_{t}\|_{L^{2}(\Omega)}^{2}.
 \end{equation}
Combining \eqref{3.8}, \eqref{3.9}, \eqref{3.10}, \eqref{3.11}, \eqref{3.13}, \eqref{3.14}, \eqref{3.15} and \eqref{3.16}, we obtain
 \begin{equation}\nonumber
\begin{split}
\frac{1}{2}\frac{d}{dt}&(\|v_{m,t}\|_{L^{2}(\Omega)}^{2}+\|d_{m,t}\|_{L^{2}(\Omega)}^{2}+\mu\|\nabla v_{m}\|_{L^{2}(\Omega_{s})}^{2})\\
&+\varepsilon \|\nabla v_{m,t}\|_{L^{2}(\Omega_{f})}^{2}+k_{1}\|\nabla d_{m,t}\|_{L^{2}(\Omega_{f})}^{2}+k_{2}\|\nabla d_{m,t}\|_{L^{2}(\Omega_{s})}^{2}\\
\leq &\frac{\varepsilon}{2}\|v_{m,t}\|_{H^{1}(\Omega_{f})}^{2}+\frac{k_{1}}{2}\|d_{m,t}\|_{H^{1}(\Omega_{f})}^{2}+\frac{1}{2}(\|f_{t}\|_{L^{2}(\Omega)}^{2}+\|g_{t}\|_{L^{2}(\Omega)}^{2})\\
&+C(T)(\|v_{m,t}\|_{L^{2}(\Omega)}^{2}+\|d_{m,t}\|_{L^{2}(\Omega)}^{2})(\|\nabla v_{m}\|_{L^{2}(\Omega_{f})}^{2}+\|\nabla d_{m}\|_{L^{2}(\Omega_{f})}^{2}+1),
\end{split}
\end{equation}
which implies
 \begin{equation}\label{3.17}
\begin{split}
\frac{d}{dt}&(\|v_{m,t}\|_{L^{2}(\Omega)}^{2}+\|d_{m,t}\|_{L^{2}(\Omega)}^{2}+\mu\|\nabla v_{m}\|_{L^{2}(\Omega_{s})}^{2})\\
&+\varepsilon \| \nabla v_{m,t}\|_{L^{2}(\Omega_{f})}^{2}+k_{1}\| \nabla d_{m,t}\|_{L^{2}(\Omega_{f})}^{2}+k_{2}\|\nabla d_{m,t}\|_{L^{2}(\Omega_{s})}^{2}\\
\leq&C(T)(\|v_{m,t}\|_{L^{2}(\Omega)}^{2}+\|d_{m,t}\|_{L^{2}(\Omega)}^{2})(\|\nabla v_{m}\|_{L^{2}(\Omega_{f})}^{2}+\|\nabla d_{m}\|_{L^{2}(\Omega_{f})}^{2}+C_{\delta})\\
&+\|f_{t}\|_{L^{2}(\Omega)}^{2}+\|g_{t}\|_{L^{2}(\Omega)}^{2},
\end{split}
\end{equation}
where $C_{\delta}=max\{\varepsilon,k_{1},k_{2},2\}$. Letting
\begin{equation*}
H(x,t)=\|v_{m,t}\|_{L^{2}(\Omega)}^{2}+\|d_{m,t}\|_{L^{2}(\Omega)}^{2}+\mu\|\nabla v_{m}\|_{L^{2}(\Omega_{s})}^{2},
\end{equation*}
then for $t=0$, and combining the estimates of \eqref{3.4} and \eqref{3.6} yields to
\begin{equation*}
H(x,0)=\|v_{m,t}(0)\|_{L^{2}(\Omega)}^{2}+\|d_{m,t}(0)\|_{L^{2}(\Omega)}^{2}+\mu\|\nabla v_{0}\|_{L^{2}(\Omega_{s})}^{2}\leq C<\infty.
\end{equation*}
Applying Gronwall's inequality to \eqref{3.17} and using assumption \eqref{1.14}, we get
\begin{equation}\label{3.18}
\begin{split}
H(x,t)\leq &e^{C(T)\int_{0}^{T}(\|\nabla v_{m}\|_{L^{2}(\Omega_{f})}^{2}+\|\nabla d_{m}\|_{L^{2}(\Omega_{f})}^{2}+C_{\delta})dt}(H(x,0)\\
&+\int_{0}^{t}(\|f_{t}(0)\|_{L^{2}(\Omega)}^{2}+\|g_{t}(0)\|_{L^{2}(\Omega)}^{2})dt\\
\leq& C(T)<\infty,\; t\in[0,T].
\end{split}
\end{equation}
Together with \eqref{3.17}, \eqref{3.18} and the energy estimates \eqref{2.4}, we obtain a prior estimates of $(v_{m,t}(t), d_{m,t}(t))$ uniformly on $m=1,2,\cdots$ that\\
\begin{equation}\label{3.19}
\|v_{m,t}\|_{L^{2}(\Omega)}^{2}+\|d_{m,t}\|_{L^{2}(\Omega)}^{2}+\|v_{m}\|_{V_{s}}^{2}+\| v_{m,t}\|_{L^{2}(0,T;V_{f})}^{2}+\| d_{m,t}\|_{L^{2}(0,T;H^{1}(\Omega))}^{2}\leq C(T).
\end{equation}\\
\textbf{Step 3: A prior estimates of  $(D_{h}v_{m}(t), D_{h}d_{m}(t))$.}\\

In this section, we consider the domain $\Omega_{f}=\mathcal{T}\times(-\tfrac{L}{2}, \tfrac{L}{2})$, $\Omega_{s}=\mathcal{T}\times((-L,-\tfrac{L}{2})\cup(\tfrac{L}{2},L))$ with the torus $\mathcal{T}=\tfrac{\mathbb{R}}{2\pi\mathcal{Z}}$, $\Gamma=\{(x,y)|x\in\mathcal{T},y=-\tfrac{L}{2}\; and\; y=\tfrac{L}{2}\}$, $d\Gamma=ds=dx$ on $x\in\Gamma$, $n=\{0, \pm 1\}$ and $\partial\Omega=\{(x,y)|x\in\mathcal{T},y=-L\; and\; y=L\}$.
On the one hand, letting $\varphi_{j}=D_{-h}D_{h}v_{m}$ in $\eqref{2.2}_{1}$, we obtain
\begin{align}\label{3.20}
\begin{split}
&(v_{m,t},D_{-h}D_{h}v_{m})_{\Omega}+\varepsilon a_{\Omega_{f}}(v_{m},D_{-h}D_{h}v_{m})+\mu a_{\Omega_{s}}(\int_{0}^{t}v_{m}d \sigma+w_{0}, D_{-h}D_{h}v_{m})\\
&+b_{\Omega_{f}}(D_{h}v_{m},v_{m},D_{h}v_{m})
-\frac{1}{2}\int_{\Gamma}(D_{h}v_{m}\cdot v_{m}) D_{h}v_{m}\cdot nd\Gamma\\
=&(d_{m} e,D_{-h}D_{h}v_{m})_{\Omega_{f}}+(f,D_{-h}D_{h}v_{m})_{\Omega},\\
\end{split}
\end{align}
where we have used the following result that
\begin{align}\nonumber
\begin{split}
&b_{\Omega_{f}}(v_{m},v_{m},D_{-h}D_{h}v_{m})-\frac{1}{2}\int_{\Gamma}(v_{m}\cdot n)(v_{m}\cdot D_{-h}D_{h}v_{m})d\Gamma\\
=&hb_{\Omega_{f}}(D_{h}v_{m},D_{h}v_{m},D_{h}v_{m})+b_{\Omega_{f}}(D_{h}v_{m},v_{m},D_{h}v_{m})+b_{\Omega_{f}}(v_{m},D_{h}v_{m},D_{h}v_{m})\\
&-\frac{1}{2}\int_{\Gamma}[h|D_{h}v_{m}|^{2}D_{h}v_{m}\cdot n+(D_{h}v_{m}\cdot v_{m}) D_{h}v_{m}\cdot n+|D_{h}v_{m}|^{2}v_{m}\cdot n]d\Gamma\\
=&b_{\Omega_{f}}(D_{h}v_{m},v_{m},D_{h}v_{m})-\frac{1}{2}\int_{\Gamma}(D_{h}v_{m}\cdot v_{m}) D_{h}v_{m}\cdot nd\Gamma.\\
\end{split}
\end{align}
Therefore, applying the property of second difference quotient $\eqref{1.11}_{1}$ and \eqref{3.20} yields
\begin{align}\label{3.21}
\begin{split}
&\frac{1}{2}\frac{d}{dt}(\|D_{h} v_{m}\|_{L^{2}(\Omega)}^{2}+\mu\|\nabla D_{h}(\int_{0}^{t}v_{m}d\sigma+w_{0m})\|_{L^{2}(\Omega_{s})}^{2})+\varepsilon \|\nabla D_{h} v_{m}\|_{L^{2}(\Omega_{f})}^{2}\\
=&-b_{\Omega_{f}}(D_{h}v_{m},v_{m},D_{h}v_{m})+\frac{1}{2}\int_{\Gamma}( D_{h}v_{m}\cdot v_{m})(D_{h}v_{m}\cdot n)d\Gamma\\
&+(D_{h}(d e),D_{h}v_{m})_{\Omega_{f}}
+(D_{h}f,D_{h}v_{m})_{\Omega},\\
=&\sum_{i=1}^{4}M_{i}.
\end{split}
\end{align}
We will estimate each term on the right-hand side of equation \eqref{3.21}.
\begin{equation}\label{3.22}
\begin{split}
M_{1}
\leq&C\|D_{h}v_{m}\|^{2}_{L^{4}(\Omega_{f})}\|\nabla v_{m}\|_{L^{2}(\Omega_{f})}\\
\leq&C\|D_{h}v_{m}\|_{H^{1}(\Omega_{f})}\|D_{h}v_{m}\|_{L^{2}(\Omega_{f})}\|\nabla v_{m}\|_{L^{2}(\Omega_{f})}\\
\leq&\frac{\varepsilon}{8}\|D_{h}v_{m}\|^{2}_{H^{1}(\Omega_{f})}+C\|D_{h}v_{m}\|^{2}_{L^{2}(\Omega_{f})}\|\nabla v_{m}\|^{2}_{L^{2}(\Omega_{f})},
\end{split}
\end{equation}
where we have used H\"{o}lder's inequality and Young's inequality. Combining the trace theorem and interpolation inequality with the 2D case, we can infer
\begin{equation}\label{3.23}
\begin{split}
M_{2}\leq&C\|D_{h}v_{m}|_{\Gamma}\|^{2}_{L^{3}(\Gamma)}\|v_{m}|_{\Gamma}\|_{L^{3}(\Gamma)}\\
\leq&C\|D_{h}v_{m}\|^{2}_{H^{2/3}(\Omega_{f})}\|v_{m}\|_{H^{2/3}(\Omega_{f})}\\
\leq&C\|D_{h}v_{m}\|^{4/3}_{H^{1}(\Omega_{f})}\|D_{h}v_{m}\|^{2/3}_{L^{2}(\Omega_{f})}\|v_{m}\|^{2/3}_{H^{1}(\Omega_{f})}\|v_{m}\|^{1/3}_{L^{2}(\Omega_{f})}\\
\leq&\frac{\varepsilon}{8}\| D_{h}v_{m}\|^{2}_{H^{1}(\Omega_{f})}+C\|D_{h}v_{m}\|^{2}_{L^{2}(\Omega_{f})}\|v_{m}\|^{2}_{H^{1}(\Omega_{f})}\|v_{m}\|_{L^{2}(\Omega_{f})},\\
\leq&\frac{\varepsilon}{8}\| D_{h}v_{m}\|^{2}_{H^{1}(\Omega_{f})}+C(T)\|D_{h}v_{m}\|^{2}_{L^{2}(\Omega_{f})}\|\nabla v_{m}\|^{2}_{L^{2}(\Omega_{f})},\\
\end{split}
\end{equation}
and
\begin{equation}\label{3.24}
M_{3}+M_{4}\leq\|D_{h}v_{m}\|^{2}_{L^{2}(\Omega)}+\frac{1}{2}(\|D_{h}d_{m}\|^{2}_{L^{2}(\Omega_{f})}+\|D_{h}f\|^{2}_{L^{2}(\Omega)}).
\end{equation}
Substituting the above estimates \eqref{3.22}, \eqref{3.23}, \eqref{3.24} into \eqref{3.21} gives
\begin{equation}\label{3.25}
\begin{split}
&\frac{1}{2}\frac{d}{dt}(\|D_{h} v_{m}\|_{L^{2}(\Omega)}^{2}+\mu\|\nabla D_{h}( \int_{0}^{t}v_{m}d\sigma+w_{0m})\|_{L^{2}(\Omega_{s})}^{2})+\varepsilon \|\nabla D_{h} v_{m}\|_{L^{2}(\Omega_{f})}^{2}\\
\leq&\frac{\varepsilon}{4}\|D_{h} v_{m}\|_{H^{1}(\Omega_{f})}^{2}+C\|D_{h}v_{m}\|_{L^{2}(\Omega)}^{2}(\| \nabla v_{m}\|_{L^{2}(\Omega_{f})}^{2}+1)\\
&+\frac{1}{2}(\|D_{h} d_{m}\|_{L^{2}(\Omega_{f})}^{2}+\|D_{h} f\|_{L^{2}(\Omega)}^{2}).
\end{split}
\end{equation}\\
On the other hand, taking $ \phi_{j}=D_{-h}D_{h}d_{m}$ in $\eqref{2.2}_{2}$, we get formally
\begin{eqnarray}\label{3.26}
\begin{split}
&(d_{m,t},D_{-h}D_{h}d_{m})_{\Omega}+k_{1} a_{\Omega_{f}}(d_{m},D_{-h}D_{h}d_{m})+ k_{2}a_{\Omega_{s}}(d_{m},D_{-h}D_{h}d_{m})\\
+&b_{\Omega_{f}}(v_{m},d_{m},D_{-h}D_{h}d_{m})
-\frac{1}{2}\int_{\Gamma} (v_{m}\cdot n)d_{m} \cdot D_{-h}D_{h}d_{m} d\Gamma=(g,D_{-h}D_{h}d_{m})_{\Omega}.\\
\end{split}
\end{eqnarray}
In fact
\begin{eqnarray}\nonumber
\begin{split}
&b_{\Omega_{f}}(v_{m},d_{m},D_{-h}D_{h}d_{m})-\frac{1}{2}\int_{\Gamma}(v_{m}\cdot n)d_{m} \cdot D_{-h}D_{h}d_{m} d\Gamma\\
=&\int_{\Omega_{f}} v^{i}_{m}\partial^{i}d^{j}_{m}D_{-h}D_{h}d^{j}_{m}dx-\frac{1}{2}\int_{\Gamma}(v^{i}_{m} d^{j}_{m})n^{i}D_{-h}D_{h}d^{j}_{m} d\Gamma\\
=&\int_{\Omega_{f}} D_{h}(v^{i}_{m}\partial^{i}d^{j}_{m})D_{h}d^{j}_{m}dx-\frac{1}{2}\int_{\Gamma}D_{h}(v^{i}_{m} d^{j}_{m})n^{i}D_{h}d^{j}_{m} d\Gamma\\
=&\int_{\Omega_{f}} hD_{h}v^{i}_{m}\partial^{i}D_{h}d^{j}_{m}D_{h}d^{j}_{m}dx+\int_{\Omega_{f}} D_{h}v^{i}_{m}\partial^{i}d_{m}^{j}D_{h}d^{j}_{m}dx+\int_{\Omega_{f}}v^{i}_{m}\partial^{i}D_{h}d^{j}_{m}D_{h}d^{j}_{m}dx\\
&-\frac{1}{2}\int_{\Gamma}(hD_{h}v^{i}_{m}D_{h}d^{j}_{m}n^{i}D_{h}d^{j}_{m}+ D_{h}v^{i}_{m}d^{j}_{m}n^{i}D_{h}d^{j}_{m}+v^{i}_{m}D_{h}d^{j}_{m}n^{i} D_{h}d^{j}_{m})d\Gamma\\
=&\int_{\Omega_{f}} D_{h}v^{i}_{m}\partial^{i}d_{m}^{j}D_{h}d^{j}_{m}dx-\frac{1}{2}\int_{\Gamma}D_{h}v^{i}_{m}d^{j}_{m}n^{i}D_{h}d^{j}_{m}d\Gamma,
\end{split}
\end{eqnarray}
hence, which together with the property $\eqref{1.11}$, \eqref{3.26} becomes
\begin{equation}\label{3.27}
\begin{split}
&\frac{1}{2}\frac{d}{dt}\|D_{h}d_{m}\|_{L^{2}(\Omega)}^{2}+k_{1} \|\nabla D_{h} d_{m}\|_{L^{2}(\Omega_{f})}^{2}+k_{2} \|\nabla D_{h} d_{m}\|_{L^{2}(\Omega_{s})}^{2}\\
=&-\int_{\Omega_{f}} D_{h}v^{i}_{m}\partial^{i}d_{m}^{j}D_{h}d^{j}_{m}dx+\frac{1}{2}\int_{\Gamma}D_{h}v^{i}_{m}d^{j}_{m}n^{i}D_{h}d^{j}_{m}d\Gamma+\int_{\Omega} D_{h}g \cdot D_{h}d_{m} dx,\\
\triangleq&\sum_{i=1}^{3}N_{i}.
\end{split}
\end{equation}
We will obtain the estimates of $N_{i}(i=1,2,3)$. For $N_{1}$, by H\"{o}lder's inequality, Young's inequality and \eqref{1.12} with $\mathrm{d}=2$, we get
\begin{equation}\label{3.28}
\begin{split}
N_{1}&=-\int_{\Omega_{f}} D_{h}v^{i}_{m}\partial^{i}d_{m}^{j}D_{h}d^{j}_{m}dx\\
\leq&C\|D_{h}v_{m}\|_{L^{4}(\Omega_{f})}\|\nabla d_{m}\|_{L^{2}(\Omega_{f})}\|D_{h}d_{m}\|_{L^{4}(\Omega_{f})}\\
\leq&C\|D_{h}v_{m}\|^{1/2}_{H^{1}(\Omega_{f})}\|D_{h}v_{m}\|^{1/2}_{L^{2}(\Omega_{f})}\|D_{h}d_{m}\|^{1/2}_{H^{1}(\Omega_{f})}\| D_{h}v_{m}\|^{1/2}_{H^{1}(\Omega_{f})}\|\nabla d_{m}\|_{L^{2}(\Omega_{f})}\\
\leq & \frac{\varepsilon}{8}\|D_{h}v_{m}\|_{H^{1}(\Omega_{f})}^{2}+\frac{k_{1}}{4}\|D_{h}d_{m}\|_{H^{1}(\Omega_{f})}^{2}+C(\|D_{h}v_{m}\|^{2}_{L^{2}(\Omega_{f})}+\|D_{h}d_{m}\|^{2}_{L^{2}(\Omega_{f})})\|\nabla d_{m}\|^{2}_{L^{2}(\Omega_{f})},\\
\end{split}
\end{equation}
As for $N_{2}$, using trace theorem and $\eqref{1.12}_{2}$, with the help of H\"{o}lder's interpolation inequality and Young's inequality, we have
\begin{equation}\label{3.29}
\begin{split}
N_{2}=&\frac{1}{2}\int_{\Gamma}D_{h}v^{i}_{m}d^{j}_{m}n^{i}D_{h}d^{j}_{m}d\Gamma\\
\leq&C\|D_{h}v_{m}|_{\Gamma}\|_{L^{3}(\Gamma)}\|D_{h}d_{m}|_{\Gamma}\|_{L^{3}(\Gamma)}\|d_{m}|_{\Gamma}\|_{L^{3}(\Gamma)}\\
\leq&C\|D_{h}v_{m}\|_{H^{2/3}(\Omega_{f})}\|D_{h}d_{m}\|_{H^{2/3}(\Omega_{f})}\|d_{m}\|_{H^{2/3}(\Omega_{f})}\\
\leq&C\|D_{h}v_{m}\|^{2/3}_{H^{1}(\Omega_{f})}\|D_{h}v_{m}\|^{1/3}_{L^{2}(\Omega_{f})}\|D_{h}d_{m}\|^{2/3}_{H^{1}(\Omega_{f})}\|D_{h}d_{m}\|^{1/3}_{L^{2}(\Omega_{f})}
\|d_{m}\|^{2/3}_{H^{1}(\Omega_{f})}\|d_{m}\|^{1/3}_{L^{2}(\Omega_{f})}\\
\leq&\frac{\varepsilon}{8}\|D_{h}v_{m}\|^{2}_{H^{1}(\Omega_{f})}+\frac{k_{1}}{4}\|D_{h}d_{m}\|^{2}_{H^{1}(\Omega_{f})}+C(T)(\|D_{h}v_{m}\|^{2}_{L^{2}(\Omega_{f})}+\|D_{h}d_{m}\|^{2}_{L^{2}(\Omega_{f})})\|\nabla d_{m}\|^{2}_{L^{2}(\Omega_{f})}.
\end{split}
\end{equation}
It is easy to derive that
\begin{equation}\label{3.30}
N_{3}\leq\frac{1}{2}(\|D_{h}g\|^{2}_{L^{2}(\Omega)}+\|D_{h}d_{m}\|^{2}_{L^{2}(\Omega)}).
\end{equation}
We plug \eqref{3.28}, \eqref{3.29}, \eqref{3.30} into \eqref{3.27} to get
\begin{equation}\label{3.31}
\begin{split}
&\frac{1}{2}\frac{d}{dt}\|D_{h} d_{m}\|_{L^{2}(\Omega)}^{2}+k_{1} \|\nabla D_{h} d_{m}\|_{L^{2}(\Omega_{f})}^{2}+k_{2} \|\nabla D_{h}d_{m}\|_{L^{2}(\Omega_{s})}^{2}\\
\leq&\frac{\varepsilon}{4}\|D_{h} v_{m}\|_{H^{1}(\Omega_{f})}^{2}+\frac{k_{1}}{2}\|D_{h} d_{m}\|_{H^{1}(\Omega_{f})}^{2}+\frac{1}{2}\|D_{h}d_{m}\|_{L^{2}(\Omega)}^{2}+\frac{1}{2}\|D_{h}g\|^{2}_{L^{2}(\Omega)}\\
&+C(T)(\|D_{h}v_{m}\|^{2}_{L^{2}(\Omega_{f})}+\|D_{h}d_{m}\|^{2}_{L^{2}(\Omega_{f})})\|\nabla d_{m}\|^{2}_{L^{2}(\Omega_{f})}.
\end{split}
\end{equation}
Collecting \eqref{3.25} and \eqref{3.31} yields
\begin{equation*}
\begin{split}
&\frac{1}{2}\frac{d}{dt}(\|D_{h} v_{m}\|_{L^{2}(\Omega)}^{2}+\mu\|\nabla D_{h}(\int_{0}^{t}v_{m}d\sigma+w_{0m})\|_{L^{2}(\Omega_{s})}^{2}+\|D_{h} d_{m}\|_{L^{2}(\Omega)}^{2})\\
&+\varepsilon \|\nabla D_{h} v_{m}\|_{L^{2}(\Omega_{f})}^{2}+k_{1} \|\nabla D_{h} d_{m}\|_{L^{2}(\Omega_{f})}^{2}+k_{2} \|\nabla D_{h} d_{m}\|_{L^{2}(\Omega_{s})}^{2}\\
\leq&\frac{\varepsilon}{2}\| D_{h}v_{m}\|^{2}_{H^{1}(\Omega_{f})}+\frac{k_{1}}{2}\| D_{h}d_{m}\|^{2}_{H^{1}(\Omega_{f})}+\|D_{h}d_{m}\|^{2}_{L^{2}(\Omega)}+C\|D_{h}v_{m}\|^{2}_{L^{2}(\Omega)}(\|\nabla v_{m}\|^{2}_{L^{2}(\Omega_{f})}+1)\\
&+C(T)(\|D_{h}v_{m}\|^{2}_{L^{2}(\Omega_{f})}+\|D_{h} d_{m}\|^{2}_{L^{2}(\Omega_{f})})\|\nabla d_{m}\|^{2}_{L^{2}(\Omega_{f})}+\frac{1}{2}(\|D_{h} f\|_{L^{2}(\Omega)}^{2}+\|D_{h} g\|_{L^{2}(\Omega)}^{2}).
\end{split}
\end{equation*}
Furthermore, we obtain
\begin{equation}\label{3.32}
\begin{split}
\frac{d}{dt}&(\|D_{h} v_{m}\|_{L^{2}(\Omega)}^{2}+\mu\|\nabla D_{h}(\int_{0}^{t}v_{m}d\sigma+w_{0})\|_{L^{2}(\Omega_{s})}^{2}+\|D_{h} d_{m}\|_{L^{2}(\Omega)}^{2})\\
&+\varepsilon\|D_{h} v_{m}\|_{H^{1}(\Omega_{f})}^{2}+k_{1}\|D_{h} d_{m}\|_{H^{1}(\Omega_{f})}^{2}+k_{2}\|D_{h} d_{m}\|_{H^{1}(\Omega_{s})}^{2}\\
\leq&C(T)(\|D_{h} v_{m}\|_{L^{2}(\Omega)}^{2}+\|D_{h} d_{m}\|_{L^{2}(\Omega)}^{2})(\|\nabla v_{m}\|_{L^{2}(\Omega_{f})}^{2}+\|\nabla d_{m}\|_{L^{2}(\Omega_{f})}^{2}+C_{\delta}).\\
&+\|D_{h} f\|_{L^{2}(\Omega)}^{2}+\|D_{h} g\|_{L^{2}(\Omega)}^{2}.
\end{split}
\end{equation}\\
By the fact $\|d_{t}\|_{L^{2}(0,T;H^{1}(\Omega))}\leq C$ and $d=d_{0}+\int_{0}^{t}d_{t}d\tau$, we obtain
\begin{equation}\label{3.33}
\|d_{m}\|_{L^{\infty}(0,T;H^{1}(\Omega))}\leq C.
\end{equation}
Similarly,
\begin{equation}\label{3.34}
\begin{split}
\|v_{m}\|_{L^{\infty}(0,T;V)}\leq C.
\end{split}
\end{equation}
Adopting Gronwall's inequality to \eqref{3.32}, combining the estimate \eqref{3.33}, the estimate \eqref{3.34} and the assumptions about $(f(x,t),g(x,t))$ stated in \eqref{1.14}, which yields
\begin{equation}\label{3.35}
\begin{split}
&\|D_{h} v_{m}\|_{L^{2}(\Omega)}^{2}+\mu\|\nabla D_{h}(\int_{0}^{t}v_{m}d\sigma+w_{0m})\|_{L^{2}(\Omega_{s})}^{2}+\|D_{h} d_{m}\|_{L^{2}(\Omega)}^{2}\\
&+\varepsilon\|D_{h} v_{m}\|_{L^{2}(0,T;H^{1}(\Omega_{f}))}^{2}+k_{1}\|D_{h} d_{m}\|_{L^{2}(0,T;H^{1}(\Omega_{f}))}^{2}+k_{2}\|D_{h} d_{m}\|_{L^{2}(0,T;H^{1}(\Omega_{s}))}^{2}\\
\leq&C(T).
\end{split}
\end{equation}
Applying the trace theorem to \eqref{3.35}, we have
\begin{equation*}
\|D_{h}(\int_{0}^{t}v_{m}d\sigma+w_{0m})\|_{H^{1/2}(\Gamma)}^{2}
+\|D_{h} v_{m}\|_{L^{2}(0,T;H^{1/2}(\Gamma))}^{2}+\|D_{h} d_{m}\|_{L^{2}(0,T;H^{1/2}(\Gamma))}^{2}\leq C(T).
\end{equation*}
Thus, we get that
\begin{equation}\nonumber
\|\int_{0}^{t}v_{m}d\sigma+w_{0}\|_{H^{3/2}(\Gamma)}^{2}
+\|v_{m}\|_{L^{2}(0,T;H^{3/2}(\Gamma))}^{2}+\| d_{m}\|_{L^{2}(0,T;H^{3/2}(\Gamma))}^{2}
\leq C(T).
\end{equation}
Thanks to
\begin{equation}\label{3.36}
\begin{cases}
-\varepsilon\Delta v_{m}+\nabla p_{m}=-v_{m,t}-v_{m}\cdot\nabla v_{m}+d_{m} e+f, &x\in\Omega_{f}, \\
\text{div}v_{m}=0, &x\in\Omega_{f}, \\
-k_{1}\Delta d_{m}=-d_{m,t}-v_{m}\cdot\nabla d_{m}+g, &x\in\Omega_{f}, \\
\mu\Delta(\int_{0}^{t}v_{m}d\sigma+w_{0m})=v_{m,t}-f,& x\in\Omega_{s}, \\
-k_{2}\Delta d_{m}=-d_{m,t}+g, &x\in\Omega_{s}. \\
 \varepsilon\frac{\partial v_{m}}{\partial n}-p_{m}n-\frac{1}{2}(v_{m}\cdot n)v_{m}=\mu\frac{\partial}{\partial n}(\int_{0}^{t}v_{m}d\sigma+w_{0m}),&(x,t)\in\Gamma\times(0,T),\\
  k_{1}\frac{\partial d_{m}}{\partial n}-\frac{1}{2} (v_{m}\cdot n)d=k_{2}\frac{\partial d_{m}}{\partial n},&(x,t)\in\Gamma\times(0,T),\\
  w_{0m}+\int_{0}^{t}v_{m} d\sigma=0,\,\frac{\partial d_{m}}{\partial N}=0, &(x,t)\in\Gamma_{out}\times(0,T).
\end{cases}
\end{equation}
Based on the regularity theorem of Stokes equation in \cite{refLadyzhenskaya} to $\eqref{3.36}_{1}$, we arrive at
\begin{equation}\label{3.37}
\begin{split}
&\|v_{m}\|_{H^{2}(\Omega_{f})}+\|p_{m}\|_{H^{1}(\Omega_{f})}\\
\leq&\|v_{m,t}\|_{L^{2}(\Omega_{f})}+\|v_{m}\cdot\nabla v_{m}\|_{L^{2}(\Omega_{f})}+\|d_{m} e\|_{L^{2}(\Omega_{f})}+\|f\|_{L^{2}(\Omega_{f})}+\|v|_{\Gamma}\|_{H^{3/2}(\Gamma)}\\
\leq&C+\|v_{m}\|_{L^{4}(\Omega_{f})}\|\nabla v_{m}\|_{L^{4}(\Omega_{f})}+\|d_{m}\|_{L^{2}(\Omega_{f})}+\|v_{m}|_{\Gamma}\|_{H^{3/2}(\Gamma)}\\
\leq&C+C\|v_{m}\|_{L^{4}(\Omega_{f})}\|\nabla v_{m}\|^{1/2}_{L^{2}(\Omega_{f})}\|\nabla v_{m}\|^{1/2}_{H^{1}(\Omega_{f})}+\|d_{m} \|_{L^{2}(\Omega_{f})}+\|v_{m}|_{\Gamma}\|_{H^{3/2}(\Gamma)}\\
\leq&C+C\|v_{m}\|^{2}_{L^{4}(\Omega_{f})}\|\nabla v_{m}\|_{L^{2}(\Omega_{f})}+\|v_{m}|_{\Gamma}\|_{H^{3/2}(\Gamma)}.\\
\end{split}
\end{equation}
Then, we deduce that
\begin{equation}\label{3.38}
\begin{split}
&\|v_{m}\|_{L^{2}(0,T;{H^{2}(\Omega_{f})})}+\|p_{m}\|_{L^{2}(0,T;{H^{1}(\Omega_{f})})}\\
\leq&C+C\int_{0}^{t}\|v_{m}\|^{4}_{L^{4}(\Omega_{f})}\|\nabla v_{m}\|^{2}_{L^{2}(\Omega_{f})}ds\\
\leq&C+C\|v_{m}\|^{4}_{L^{\infty}(0,T;L^{4}(\Omega_{f}))}\|\nabla v_{m}\|^{2}_{L^{2}(0,T;L^{2}(\Omega_{f}))}\\
\leq&C+C(\|v_{0}\|^{4}_{L^{4}(\Omega_{f})}+\|v_{m,t}\|^{4}_{L^{2}(0,T;L^{4}(\Omega_{f}))})\\
\leq&C(T).
\end{split}
\end{equation}
Using a similar argument to $\eqref{3.36}_{3}$ and $\eqref{3.36}_{5}$, we get that
\begin{equation}\label{3.39}
\|d_{m}\|_{L^{2}(0,T;{H^{2}(\Omega)})}\leq C(T).
\end{equation}
By the regularity theorem of elliptic equation in \cite{refLions1} to $\eqref{3.36}_{4}$, we have
\begin{equation}\label{3.40}
\|\int_{0}^{t}v_{m}d\sigma+w_{0m}\|_{H^{2}(\Omega_{s})}\leq\|v_{m,t}\|_{L^{2}(\Omega_{s})}+\|f\|_{L^{2}(\Omega_{s})}+\|\int_{0}^{t}v_{m}d\sigma+w_{0m}\|_{H^{3/2}(\Gamma)}\leq C(T).
\end{equation}\\
Now, we get the proof of \eqref{1.16} when passing the limit $m\rightarrow \infty$ in \eqref{3.19}, \eqref{3.38}, \eqref{3.39} and \eqref{3.40}.
Furthermore, we can establish the more high regularity when the smooth initial value satisfies the more high compatibility conditions.
Theorem \ref{thm 1.2.} has been completed.

\section{\bf Acknowledgements}
This work was partially supported by National Natural Science Foundation of China (No.12171111).

\section{\bf Conflict of Interest}
This work does not have any conflicts of interest.

 \small{
}

\begin{thebibliography}{99}


\bibitem{refBreit} D. Breit, S. Schwarzacher, Compressible fluids interacting with a linear-elastic shell, Arch. Ration. Mech. Anal., 228(2018) 495-562.

\bibitem{refBoffi} D. Boffi, L. Gastaldi, On the existence and the uniqueness of the solution to a fluid-structure interaction problem, J. Differential
Equations, 279(2021) 136-161.

\bibitem{refBoyer} F. Boyer, P. Fabrie, Mathematical Tools for the Study of the Incompressible Navier-Stokes
Equations and Related Models. Applied Mathematical Sciences, vol. 183, Springer, New-York, 2013.

\bibitem{refBae} H. Bae,: On the free boundary value problem of the Navier-Stokes equations. Discrete Contin. Dyn. Syst., 29(2011) 769-801.

\bibitem{refHB} H. Beir\~{a}o da Veiga, On the existence of strong solutions to a coupled fluid-structure evolution problem, J. Math. Fluid Mech., 6(1)(2004) 21-52.

\bibitem{refBarbu1} V. Barbu, Z. Gruji\'{c}, I. Lasiecka, A. Tuffaha, Existence of the energy-level weak solutions for a nonlinear fluid-structure interaction model, in: Fluids and Waves, in: Contemp. Math., vol. 440, Amer. Math. Soc., Providence, RI, 2007, pp. 55-82.

\bibitem{refBarbu2} V. Barbu, Z. Gruji\'{c}, I. Lasiecka, A. Tuffaha, Smoothness of weak solutions to a nonlinear fluid-structure interaction model, Indiana Univ. Math. J., 57(3)(2008) 1173-1207.

\bibitem{refBodnar} T. Bodn\'{a}r, G. P. Galdi, and \v{S}. Ne\v{c}asov\'{a} (eds.), Fluid-structure interaction and biomedical applications. Advances in Mathematical Fluid Mechanics, Birkh\"{a}user/Springer, Basel, 2014 Zbl 1300.76003 MR 3223031

\bibitem{refChambolle} A. Chambolle, B. Desjardins, M.J. Esteban, C. Grandmont, Existence of weak solutions for the
unsteady interaction of a viscous fluid with an elastic plate, J. Math. Fluid Mech., 7(3)(2005) 368-404.

\bibitem{refCoutand1} D. Coutand and S. Shkoller, Motion of an elastic solid inside an incompressible viscous fluid, Arch. Ration. Mech. Anal., 176(1)(2005) 25-102.

\bibitem{refCoutand2} D. Coutand and S. Shkoller, The interaction between quasilinear elastodynamics and the Navier-Stokes equations, Arch. Ration.
    Mech. Anal., 179(3)(2006) 303-352.

\bibitem{refC} H. Cohen and S.I. Rubinow, Some mathematical topics in biology, Proc. Symp. on System Theory, Polytechnic Press, New York (1965), 321-337.

\bibitem{refCrosetto} P. Crosetto, P. Reymond, S. Deparis, D. Kontaxakis, N. Stergiopulos, A.Quarteroni, Fluid-structure interaction simulation of aortic blood flow. Comput Fluids, 43(1)(2011) 46-57.

\bibitem{refDu} Q. Du, M.D. Gunzburger, L.S. Hou, J. Lee, Analysis of a linear fluid-structure interaction problem, Discrete Contin.
Dyn. Syst., 9(3)(2003)633-650.

\bibitem{refGrandmont} C. Grandmont, Existence of weak solutions for the unsteady interaction of a viscous fluid with an elastic plate,
 SIAM J. Math. Anal., 40(2)(2008) 716-737.

\bibitem{refDowell} E.H. Dowell, A modern course in aeroelasticity. Enl. edn., Solid Mechanics Appl 217,
Springer, Cham, 2015 Zbl 1297.74001 MR 3306893

\bibitem{refKukavica1} I. Kukavica, A. Tuffaha, M. Ziane, Strong solutions to a nonlinear fluid structure interaction system, J. Differ. Equations, 247(5)(2009) 1452-1478.

\bibitem{refKukavica} I. Kukavica, A. Tuffaha, M. Ziane, Strong solutions for a fluid structure interaction system, Adv. Differential Equations, 15(3-4)(2010) 231-254.

 \bibitem{refKuttler} U. K\"{u}ttler, W.A. Wall, The dilemma of domain decomposition approaches in fluid-structure interactions with fully enclosed incompressible fluids. In: Domain decomposition methods in science and engineering XVII. Lect Notes Comput Sci Eng, vol. 60. Berlin: Springer, 2008. p. 575-582

\bibitem{refLengeler1} D. Lengeler, Weak solutions for an incompressible, generalized  Newtonian  fluid  interacting with a linearly elastic  Koiter  type shell, SIAM J. Math. Anal., 46(2014) 2614-2649

\bibitem{refLadyzhenskaya} O.A. Ladyzhenskaya, Solution "in the large" to the boundary value problem for the Navier-Stokes equations in two space variables. Sov. Phys. Dokl., 3(1958), 1128-1131. Translation from Dokl. Akad. Nauk SSSR 123 (1958) 427-429.

\bibitem{refLasiecka} I. Lasiecka, J.L. Lions, R. Triggiani, Nonhomogeneous boundary value problems for second order hyperbolic operators, J. Math. Pures Appl., 65(2)(1986) 149-192.

\bibitem{refLequeurre1} J. Lequeurre, Existence of strong solutions to a fluid-structure system, SIAM J. Math. Anal., 43(1)(2011) 389-410.

\bibitem{refLions} J. L. Lions, Quelques M\'{e}thodes de R\'{e}solution des Problmes aux Limites non Lin\'{e}aires. Paris: Dunod, 1969.

\bibitem{refLions1} J.L. Lions, E. Magenes, Non-Homogeneous Boundary Value Problems and Applications I, Springer-Verlag, Berlin, New York, 1972, pp. 152-155.

\bibitem{refShen} L. Shen, S. Wang, Existence of global weak solutions for the high frequency and small displacement oscillation fluid-structure interactions systems, Math Meth Appl Sci., (2020) 1-11.

\bibitem{refShen1} L. Shen, S. Wang, Note on the global wellposedness of two-dimensional incompressible Navier-Stokes fluid structure interaction model. https://doi.org/10.1016/j.aml.2022.108564.

\bibitem{refLx}  L.X. Zhang, Y.K. Guo, H.M. Zhang, Fully coupled flow-induced vibration of structures under small deformation with GMRES method. Appl. Math. Mech.-Engl. Ed., 31(2010) 87-96.

\end{thebibliography}
\end{document}